\documentclass[twocolumn,draft]{autart}
\usepackage[pdftex,bookmarks,colorlinks,citecolor=black linkcolor=blue]{hyperref}
\usepackage[pdftex,final]{graphicx}
\usepackage{amsmath,amssymb,amsfonts,dsfont}
\usepackage{mathrsfs} 

\begin{document}
\begin{frontmatter}
\runtitle{Minimax estimates for linear DAEs}
\title{Minimax state estimation for linear continuous differential-algebraic equations}
\author{Sergiy M. Zhuk}
\address{Department of System Analysis and Decision Making Theory, Taras Shevchenko National University of Kyiv, Ukraine} 
\ead{beetle@unicyb.kiev.ua}
%\date{10 January 2011}
%\maketitle

\begin{abstract}
%\noindent{\bf Abstract}\quad 
This paper describes a minimax state estimation approach for linear Differential-Algebraic Equations (DAE) with uncertain parameters. The approach addresses continuous-time DAE with non-stationary rectangular matrices and uncertain bounded deterministic input. An observation's noise is supposed to be random with zero mean and unknown bounded correlation function. 
%The state's observations are supposed to be incomplete and noisy; the noise is modelled by a random process with zero mean and unknown bounded correlation function. \\
Main results are a Generalized Kalman Duality (GKD) principle and sub-optimal minimax state estimation algorithm. GKD is derived by means of Young-Fenhel duality theorem. GKD proves that the minimax estimate coincides with a solution to a Dual Control Problem (DCP) with DAE constraints. The latter is ill-posed and, therefore, the DCP is solved by means of Tikhonov regularization approach resulting a sub-optimal state estimation algorithm in the form of filter. We illustrate the approach by an synthetic example and we discuss connections with impulse-observability.  
\end{abstract}
%  and \emph{index of
%\newline\textbf{Keywords: }

% \begin{AMS}
% 34K32, %implicit ODE
% 49N10, %linear-quadratic control problems
% 49N15, %opt control duality problems
% 49N30, %problems with incomplete info
% 49N45x %optimal control, inverse problems
% \end{AMS}
\end{frontmatter}
% \pagestyle{myheadings}
%  \thispagestyle{plain}
%  \markboth{S. Zhuk}{Minimax estimates for DAE}
%TODO: Euler-Lagrange system->Hamilton system
%rivn dla filtry vypysaty jak stoh rn-nja (u dyferentsialah, Leibnits notation)!
% main text
\section{Introduction}
This paper presents a generalization of the minimax state estimation approach to linear Differential-Algebraic Equations (DAE) in the form
\begin{equation}
  \label{eq:state}
  \dfrac {d(Fx)}{dt}=C(t)x(t)+f(t), \quad Fx(t_0)=x_0 
\end{equation}
where $F\in\mathbb R^{m\times n}$ is a rectangular $m\times n$-matrix and $t\mapsto C(t)\in\mathbb R^{m\times n}$ is a continuous matrix-valued function. The research presented here may be thought of as a continuation of the paper~\cite{Zhuk2009c}, where the case of discrete time DAEs with time-depending coefficients was investigated. 
%2REWRITE The author's interest in developing a state estimation theory for systems of type~(\ref{eq:state}) is motivated, in particular, by the following computational problem coming from geophysics. The dimension of the state space in geophysical applications\cite{MalletZhuk2011} can be of order $N=10^6$. Therefore, state estimation algorithms based on the Kalman filtering approach can not\footnote{It is impossible to propagate a Riccati gain due to the computational burden.} be applied without a proper dimension reduction\cite{wu08comparison}. A reasonable and mathematically justified reduction procedure could be build by projecting the state space of the model onto some subspace $\mathcal P$ of dimension $n\ll N$. In particular, the subspace $\mathcal P$ can be computed by means of a robust Principal Component Analysis \cite{Somth on PCA}, so that the projected state represents a ``most significant part'' of the state vector. The dynamics of the resulting projected state is governed by a DAE with $F\in\mathbb R^{N\times n}$ and state space dimension equal to $n$. As $n\ll N$ one can try to estimate the state of the latter DAE (which is a projected state of the original system at the same time). Having the estimate of the projected state one can restore the ``most significant part'' of the state vector. \\
%We will illustrate the reduction procedure by an example (see Section~\ref{ss:numex}). Further details along with application to an operational air quality model could be found in \cite{MalletZhuk2011}. \\
We stress that the DAE with $F\in\mathbb R^{m\times n}$ is non-causal (the matrix pencil $F-\lambda C(t)$ is singular\cite{Gantmacher1960}) if $m\ne n$. Also the coefficient $C(t)$ depends on time. Therefore the state estimation problem for DAE in the form~\eqref{eq:state} can not be directly addressed by parameter estimation methods (see for instance \cite{Gerdin2007}, \cite{Darouach1997} and citations there), based on the transformation of the regular matrix pencil $F-\lambda C$ ($\mathrm{det}\;(F-\lambda C)\not\equiv0$) to the Weierstrass canonical form \cite{Gantmacher1960}. As it was mentioned in \cite{Gerdin2007}, the latter transformation allows to convert DAE (with regular pencil) into Ordinary Differential Equation (ODE), provided the unknown input $f$ is smooth enough and $C(t)\equiv\,\mathrm{const}$. 
On the other hand, in applications $f$ is often modelled as a realization of some random process or as a measurable squared-integrable function with bounded $\mathbb L_2$-norm. One way to go is to take Sobolev-Shvartz derivative of $f$, allowing the state $x(t)$ of DAE to be discontinuous function. If the latter is not acceptable, it is natural to ask if it is possible to derive a state estimation algorithm for DAE (in the form~(\ref{eq:state})) avoiding the differentiation of the unknown input $f$. More generally, is it possible to derive a state estimation algorithm for DAE in the form~(\ref{eq:state}) with measurable $f$ without transforming the pencil $F-\lambda C$ into a canonical form? The same question arises if $C(t)$ is not constant as in this case it may be impossible (see \cite{Campbell1987}) to transform DAE to ODE even if the pencil $F-\lambda C(t)$ is regular for all $t$. 
%Related question is connected with numerical properties of a state estimation algorithm for DAE~\eqref{eq:state}. As the solution $x$ of (\ref{eq:state}) may depend (explicitly or not) on derivatives of the input $f$ it follows that, in general case, the initial value problem for~(\ref{eq:state}) is not well-posed \cite{Hanke1989} in the space $\mathbb L_2(t_0,T,\mathbb R^n)$. The latter could introduce difficulties during the numerical integration of the ``estimation equation''. In this situation one can think about regularization of~(\ref{eq:state}) so that it becomes well-posed. However, the connection between the estimate obtained for the regularized DAE and the solution to the original ill-posed DAE should be clarified. Another workaround is to construct a sub-optimal estimation algorithm for the ill-posed DAE and to prove that it converges to the optimal algorithm. \\
In this paper we give a positive answer to this question for the following state estimation problem: given observations $y(t)$, $t\in[t_0,T]$ of $x(t)$, to reconstruct $Fx(T)$, provided $x$ is a weak solution to \eqref{eq:state}. 
% and the matrix pencil $F-\lambda C$ is not necessary regular. 
We note, that many authors (see \cite{Gerdin2007}, \cite{Darouach1997} and citations there) assume the state vector $x(t)$ of~(\ref{eq:state}) to be a differentiable (in the classical sense) function. In contrast, we only assume that $t\mapsto Fx(t)$ is an absolutely continuous function. In this setting $x(T)$ is not necessary well defined. Hence, it makes sense to estimate $Fx(T)$ only.  
%Our assumption on $x$ allows to resolve two parameter estimation problems (estimation of the state $z(t)$ and estimation of the input $f(t)$) by means of the same estimation algorithm, that is to consider the input $f(t)$ as a part of the generalized system's state $x(t) = \bigl(
% \begin{smallmatrix}
%   z(t)\\f(t)
% \end{smallmatrix}\bigr)$. The strategy of considering model's parameters (for instance input) as a part the system's state is known in the mathematical systems theory as a behavioral approach \cite{Ilchmann2005}.
In what follows, we assume that $f$ is an unknown squared-integrable function, which belongs to a given bounded set $\mathscr G$. We will also assume that observations $y(t)$ may be incomplete and noisy, that is $y(t)=H(t)x(t)+\eta(t)$, where $\eta$ is a realization of a random process with zero mean and unknown but bounded correlation function.  %$(t,s)\mapsto \mathcal R_\eta(t,s):=E\eta(t)\eta^T(s)$ belongs to a given bounded set $W$. 
Following \cite{Nakonechnii1978} we will be looking for the minimax estimate $\widehat{\ell(x)}$ of a linear function
\footnote{Note that in order to reconstruct $Fx(T)$ it is enough to reconstruct a linear function $\ell(x):=\langle\ell,Fx(T)\rangle$ for any $\ell\in\mathbb R^m$. 
Having the estimate of $\ell(x)$ for any $\ell\in\mathbb R^m$, one can set $\ell:=e_i=(0,\dots,1,\dots0)^T$ in order to reconstruct $i$-th component of $Fx(T)$.} $\ell(x):=\langle\ell,Fx(T)\rangle$ among all linear functions of observations $u(y)$. Main notions of deterministic minimax state estimation approach \cite{Tempo1985,Chernousko1994}, \cite{Kurzhanski1997} are reachability set, minimax estimate and worst-case error. By definition, reachability set contains all states of the model which are consistent with observed data and uncertainty description. Given a point $P$ within the reachability set one defines a worst-case error as the maximal distance between $P$ and other points of the reachability set. Then the minimax estimate of the state is defined as a point minimizing the worst-case error (a Tchebysheff center of the reachability set). In this paper we deal with random noise in observations. This prevents us from describing the reachability set. Instead, we derive a dynamic mean-squared minimax estimate minimizing mean-squared worst-case error. 
%By definition\cite{Nakonechnii1978}, the minimax estimate $\widehat{\ell(x)}$ minimizes a worst-case estimation error $\sigma(T,\ell,u)$, where $\sigma$ is a maximal distance $\rho$ between $u(y)$ and $\ell(x)$ with respect to uncertain parameters. Namely, $\rho$ can be defined as $E[\ell(x)- u(y)]^2$. Since $\rho$ depends (through DAE~\eqref{eq:state} and equation for $y$) on realizations of $f$, $\mathcal R_\eta$ and non-trivial solutions of homogeneous DAE~\eqref{eq:state}, it follows that $
%\sigma(T,\ell,u):=\sup_{f\in \mathscr G,\mathcal R_\eta\in W,x\text{ solves~\eqref{eq:state}}}
%\{E[\ell(x)-u(y)]^2\}
%$. Now, by definition $\widehat{\ell(x)}=\hat u(y)$, where $\sigma(T,\ell,\hat u)=\inf_u\sigma(T,\ell,u):=\hat\sigma(T,\ell)$ is referred to as a minimax error. \\
%The minimax estimate for \eqref{eq:state} was constructed in \cite{Nakonechnii1978}, provided $F=I$. The case of deterministic $\eta$ and $F=I$ was addressed in \cite{Bertsekas1971,Kurzhanski1997}, \cite{Chernousko1994} where the minimax estimation of $x(T)$ is shown to be the Tchebysheff center of the ODEs reachable set, consistent with observations and uncertainty description. % Here, we present a minimax estimation for a class of linear non-causal DAE: $F\in\mathbb R^{m\times n}$ and $t\mapsto C(t)\in\mathbb R^{m\times n}$ is continuous on $[t_0,T]$.
%The same result can be proved by the method of this paper for $F\ne I$, provided $\eta$ is deterministic and $(\eta,f)$ belongs to an ellipsoid. \\

The contributions of this paper are a Generalized Kalman Duality (GKD) and sub-optimal minimax estimation algorithm, both for DAE in the form~\eqref{eq:state}. As it was previously noted in \cite{Gerdin2007} the need to differentiate an unknown input 
%\footnote{The unknown input is typically modelled as a realization of a random process (stochastic filtering framework) or as an unknown $L_2$-function (deterministic filtering framework). 
%In the both cases the input is not a differentiable function!} 
posed a problem of mathematical justification of the filtering framework based on DAE with classical derivatives. In \cite{Gerdin2007} the authors propose a solution, provided $\mathrm{det}(F-\lambda C)\ne 0$ for any $\lambda$ \cite{Gerdin2007}. Here we apply GKD in order to justify the minimax filtering framework for the case of DAEs (\ref{eq:state}) with any rectangular $F$ and time-varying $C(t)$. We do not use the theory of matrix pencils so that the condition of differentiability of the unknown input $f$ in (\ref{eq:state}) is not necessary for our derivations. 
%At the same time the matrices $F$ and $C(t)$ are supposed to be rectangular unlike \cite{Gerdin2007}, \cite{Darouach1997}. 
%In fact, GKD gives a tool for solving the Minimax Estimation Problem (MEP) for abstract linear ill-posed equations \cite{Zhuk2009d}. In the case of DAE, the GKD states that the MEP is equal to a so-called dual control problem, provided $\ell$ belongs to so called minimax observable subspace $\mathcal L(T)$. 
%We stress that the minimax error is infinite if $\ell\not\in\mathcal L(T)$. In the latter case the estimation error $E[\ell(x)-u(y)]^2$ varies in $[0,+\infty]$ for any linear estimation $u(y)$ of $\ell(x)$. Efficient description of $\mathcal L(T)$ can be done for continuous DAE (Remark~\ref{r:1}) with constant coefficient and discrete time DAE with variable coefficients \cite{Zhuk2009c}. \\
%Some aspects of classical observability for DAE were considered in \cite{Campbell1991,Frankowska1990}. \\
Applying GKD we arrive to the Dual Control Problem (DCP) with DAE constraint, which has a unique absolutely continuous solution, provided $\ell$ belongs to the minimax observable subspace $\mathcal L(T)$. Otherwise, the solution of DCP is represented in terms of the impulsive control. In this sense the minimax observable subspace generalizes impulse observability condition (see \cite{Gerdin2007}) to the case of DAE with rectangular time-varying coefficients. The cost function of DCP describes the mean-squared worst-case error and its minimizer represents a minimax estimate. However, Pontryagin Maximum Principle (PMP) can not be applied directly to solve DCP: a straightforward application of the PMP to the dual problem could reduce the minimax observable subspace to the trivial case $\mathcal L(T)=\{0\}$ (see example in Subsection \ref{ss:optcond}). In order to preserve the structure of $\mathcal L(T)$ we apply Tikhonov regularization approach. 
%This is a consequence of the following fact: the range of the linear mapping induced by DAE is not closed\footnote{The latter corresponds to the DAE with non-zero index as the solution depends on the derivatives of the right-hand part.}\cite{Zhuk2007}. This observation suggests to apply Tikhonov method in order to find the minimax estimate. 
As a result (Proposition~\ref{p:3e}) we represent a sub-optimal minimax state estimation algorithm as a unique solution of a well-posed Euler-Lagrange system with a small parameter. This solution converges to the minimax estimate. We represent the sub-optimal estimate in the classical sequential form: as a solution to a Cauchy problem for a linear stochastic ODE, driven by a realization of observations $y(t)$, $t\in[t_0,T]$. We recall that $y(t)$ is perturbed by a ``random noise'', which can be a realization of any random process (not necessary Gaussian as in \cite{Gerdin2007}) with zero mean and unknown but bounded correlation function. \\
This paper is organized as follows. At the beginning of section~\ref{s:obs} we describe the formal problem statement and introduce definitions of the minimax mean-squared estimates and errors. The rest of this section consists of two subsections. Subsection~\ref{s:gkd} presents the GKD (Theorem~\ref{t:1}). In subsection~\ref{ss:optcond} we discuss optimality conditions and derive regularization scheme (Proposition~\ref{p:3e}) along with the representation of the sub-optimal minimax estimate in the sequential form (Corollary~\ref{c:eme}). Also we present an example. Section~\ref{sec:Conclusion} contains conclusion. Appendix contains proofs of technical statements.\\ %used in the proof of Proposition~\ref{p:3e}.\\ 
\emph{Notation}: 
$E\eta$ denotes the mean of the random element $\eta$; $\mathrm{int}\, G$ denotes the interior of $G$;  
%$f$ or $f$ denotes a function (as an element of some functional space); $f(t)$ denotes a value of function $f$ at $t$; 
%$f_t$ denotes the derivative of $t\mapsto f(t)$;
%$f^*$ denotes Young-Fenhel transformation of $f$, 
$\mathbb R^n$ denotes the $n$-dimensional Euclidean space;
$\mathbb L_2(t_0,T,\mathbb R^m)$ denotes a space of square-integrable functions with values in $\mathbb R^m$ (in what follows we will often write $\mathbb L_2$ referring $\mathbb L_2(t_0,T,\mathbb R^k)$ where the dimension $k$ will be defined by the context); 
$\mathbb H_1(t_0,T,\mathbb R^m)$ denotes a space of absolutely continuous functions with $\mathbb L_2$-derivative and values in $\mathbb R^m$; 
the prime $'$ denotes the operation of taking the adjoint: $L'$ denotes adjoint operator, $F'$ denotes the transposed matrix; 
$c(G,\cdot)$ denotes the support function of a set $G$;   %$\delta(G,\cdot)$ denotes the indicator\footnote{$\delta(G,f)=0$ if $f\in G$ and $+\infty$
%otherwise.} of $G$, %$\mathrm{dom} f=\{x:f(x)<\infty\}$; 
$\langle\cdot,\cdot\rangle$ denotes the inner product in a Hilbert space $\mathcal H$, $\|x\|^2_{\mathcal H}:=\langle x,x\rangle$, $\|\cdot\|$ denotes norm in $\mathbb R^n$; 
$S>0$ means $\langle Sx,x\rangle>0$ for all $x$; 
$F^+$ denotes the pseudoinverse matrix; $Q^\frac 12$ denotes the square-root of the symmetric non-negative matrix $Q$, $I_n$ denotes $n\times n$-identity matrix, $0_{n\times m}$ denotes $n\times m$-zero matrix, $I_0:=0$; $\mathrm{tr}$ denotes the trace of the matrix.  
\section{Linear minimax estimation for DAE}\label{s:obs}       
%We proceed with the formal problem statement. 
Consider a pair of systems% of linear differential-algebraic equations
\begin{align}
      & \dfrac{d(Fx)}{dt}=C(t)x(t)+f(t)\,,\quad Fx(t_0)=x_0\,,\label{eq:dae}\\
      & y(t)=H(t)x(t)+\eta(t)\,,    \label{eq:dae:obs}
\end{align}
where $x(t)\in\mathbb R^n$, $f(t)\in\mathbb R^m$, $y(t)\in\mathbb R^p$,
$\eta(t)\in\mathbb R^p$ represent the state, input, observation and
observation's noise respectively. As above, we assume $F\in\mathbb R^{m\times n}$, $f\in\mathbb L_2(t_0,T,\mathbb R^m)$, and $C(t)$ and $H(t)$ are continuous\footnote{Slightly modifying the proofs we can allow $C(t)$ and $H(t)$ to be just measurable.} matrix-valued functions of $t$ on $[t_0,T]$, $t_0,T\in\mathbb R$. Now let us describe our assumptions on uncertain $x_0,f,\eta$. Let $\eta$ be a realization of a random process such that $E\eta(t)=0$ on $[t_0,T]$ and $\mathcal R_\eta(t,s):=E\Psi(t)\Psi'(s)$ is bounded:
\begin{equation}
  \label{eq:eta_bounds}
\mathcal R_\eta\in W=\{\mathcal R_\eta:\int_{t_0}^T\mathrm{tr}\,( R(t)\mathcal R_\eta(t,t))dt\le 1\}
%\eta\in W=\{\eta:E\int_{t_0}^T\langle R(t)\eta(t),\eta(t)\rangle\le 1\}
\end{equation}
We note that the assumption $E\eta(t)=0$ is not restrictive. If  $E\eta(t)=\overline\eta(t)\ne 0$ where $\overline\eta (t)$ is a known measurable function then we can consider new measurements $\overline y(t):=y(t)-\overline\eta(t)$ and new noise $\eta(t)-\overline\eta(t)$. \\
The initial condition $x_0$ and input $f$ are supposed to belong to the following set
\begin{equation}
  \label{eq:G}
  \mathscr G:=\{(x_0,f):\langle Q_0x_0,x_0\rangle+\int_{t_0}^T\langle Q(t)f,f\rangle dt\le 1\},
\end{equation}
where $Q_0,Q(t)\in\mathbb R^{m\times m}$, $Q_0=Q_0'>0$, $Q=Q'>0$, $R(t)\in\mathbb R^{p\times p}$, $R'=R>0$; $Q(t)$, $R(t)$, $R^{-1}(t)$ and $Q^{-1}(t)$ are continuous functions of $t$ on $[t_0,T]$. 
\begin{defn}\label{d:1}
Given $T<+\infty$, $u\in\mathbb L_2(t_0,T,\mathbb R^p)$ and $\ell\in\mathbb R^m$ define a mean-squared worst-case estimation error 
\begin{equation*}
  \begin{split}
    &\sigma(T,\ell,u):=\sup_{x, (x_0,f)\in \mathscr G, \mathcal R_\eta\in W}
    E[\langle\ell, Fx(T)\rangle-u(y)]^2%:(Fx)_t=Cx+f, Fx(t_0)=0 \}
  \end{split}
\end{equation*}
A function $\hat u(y)=\int_{t_0}^{T}\langle\hat u(t),y(t)\rangle dt$ is called an a priori minimax mean-squared estimate in the direction $\ell$ ($\ell$-estimate) if $\inf_u\sigma(T,\ell,u)=\sigma(T,\ell,\hat u)$. The number
$\hat\sigma(T,\ell)=\sigma(T,\ell,\hat u)$ is called a minimax mean-squared a priori error in the direction $\ell$ at time-instant $T$ ($\ell$-error). 
The set $\mathcal L(T)=\{\ell\in\mathbb R^m:\hat\sigma(T,\ell)<+\infty\}$ is called a minimax observable subspace. 
\end{defn}
\subsection{Generalized Kalman Duality Principle}\label{s:gkd}
  Definition~\ref{d:1} % generalizes the notion of the linear minimax a priori mean-squared estimation, introduced in \cite{Nakonechnii1978}. It
  reflects the procedure of deriving the minimax estimation. The first step is, given $\ell$ and $u$ to calculate the worst-case error $\sigma(T,\ell,u)$ by means of the suitable duality concept. 
%, that is to formulate a dual control problem. The second step is to find the estimate $\hat u$ minimizing the worst-case error $\sigma(T,\ell,u)$, that is to solve the dual problem. The first step is given by the following theorem, which generalizes the celebrated Kalman duality principle on DAEs. %~\cite{Brammer1989} to DAE. 
\begin{thm}[Generalized Kalman duality]\label{t:1}
Take $\ell\in\mathbb R^m$. The $\ell$-error $\hat\sigma(T,\ell)$ is finite iff  
\begin{equation}
    \label{eq:zul}
    \dfrac{d(F'z)}{dt}=-C'(t)z(t)+H'(t)u(t),\quad F'z(T)=F'\ell
  \end{equation}
for some $z\in\mathbb L_2(0,T,\mathbb R^m)$ 
and $u\in\mathbb L_2(0,T,\mathbb R^p)$.\\ 
If $\hat\sigma(T,\ell)<+\infty$ then 
\begin{equation}
  \label{eq:umin}
  \begin{split}
% \sigma(T,\ell,u)=&
\sigma&(T,\ell,u)=
\min_{v,d}\{\| \tilde Q_0^{-\frac12}(z(t_0)-v(t_0))-Q_0^{-\frac12}d\|^2\\
%\min_{v,d}\{\| Q_0^{-\frac12}(({F'}^+F'(z(t_0)-v(t_0))-d)\|^2\\
&+\int_{t_0}^T\|Q^{-\frac 12}(z-v)\|^2 dt\}+
%\int_{t_0}^T
\int_{t_0}^T\| R^{-\frac12}u\|^2dt\\
&=\| \tilde Q_0^{-\frac12}(z(t_0)-\tilde v(t_0))-Q_0^{-\frac12}\tilde d\|^2\\
&+\int_{t_0}^T\|Q^{-\frac 12}(z-\tilde v)\|^2 +\| R^{-\frac12}u\|^2dt
%\quad F'd=0, v\text{obeys \eqref{eq:zul} with }u=0,\ell=0  
  \end{split}
 \end{equation}
where $\tilde Q_0^{-\frac12}=Q_0^{-\frac12}{F'}^+F'$, $\min$ in~\eqref{eq:umin} is taken over all $d$ such that $F'd=0$ and all $v$ verifying \eqref{eq:zul} with $u=0$ and $\ell=0$, and $\min$ is attained at $\tilde v,\tilde d$. 
\end{thm}
% \begin{lem}\label{l:1}
% Let $\ell\in\mathcal L(T)$. Then~\eqref{eq:umin} is equal to
% \begin{equation}
%   \label{eq:dcp}
%   N(z,u)
% :=\int_{t_0}^T\langle Q^{-1}z,z\rangle dt
% +\int_{t_0}^T\langle R^{-1}u,u\rangle dt
% \to\min_{(z,u)\text{ verify \eqref{eq:zul}}}
% \end{equation}
% \end{lem}
% The proof of Lemma~\ref{l:1} is easy and is, therefore, omitted. 
\begin{rem}\label{r:1}
  An obvious corollary of Theorem~\ref{t:1} is an expression for the minimax observable subspace
  \begin{equation}
    \label{eq:LT}
    \mathcal L(T)=\{\ell\in\mathbb R^m:\text{ \eqref{eq:zul} holds for some }z,\,u\}
  \end{equation} 
In the case of stationary $C(t)$ and $H(t)$ the minimax observable subspace may be calculated explicitly, using the canonical Kronecker form \cite{Gantmacher1960}.   
\end{rem}
\begin{pf}%{of Proposition~\ref{p:1}}
Take $\ell\in\mathbb R^m$ and $u\in\mathbb L_2(t_0,T,\mathbb R^p)$. Using $E\eta(t)=0$ we compute %$E[\ell(x)-u(y)]^2$. It is easy to see using $E\eta(t)=0$ that
\begin{equation}
  \label{eq:error}
  \begin{split}
    E[\ell&(x)-u(y)]^2 = E[\int_{t_0}^T\langle u,\eta\rangle dt]^2(:=\gamma^2)\\
&+\bigl(\langle\ell,Fx(T)\rangle-\int_{t_0}^T\langle H'u,x\rangle dt\bigr)^2(:=\mu^2)
  \end{split}
\end{equation}
Let us transform $\mu^2$. There exists $w\in\mathbb L_2(t_0,T,\mathbb R^m)$ such that: (W) $F'w\in\mathbb H^1(t_0,T,\mathbb R^n)$ and $F'w(T)=F'\ell$. Noting that \cite{Albert1972} $F=FF^+F$ we have
\begin{equation*}
%  \label{eq:lFx}
  \langle\ell,Fx(T)\rangle=\langle F'\ell,F^+Fx(T)\rangle = 
\langle F'w(T),F^+Fx(T)\rangle\,.
\end{equation*}
The latter equality, an integration-by-parts formula 
\begin{equation}
  \label{eq:ibp}
  \begin{split}
\langle F'w(T)&,F^+Fx(T)\rangle-\langle F'w(t_0),F^+Fx(t_0)\rangle=\\
 &= \int_{t_0}^T \langle\dfrac{d(Fx)}{dt},w\rangle+\langle\dfrac{d(F'w)}{dt},x\rangle dt\,. 
  \end{split}
\end{equation}
(proved in~\cite{Zhuk2007} for $Fx\in\mathbb H^1(t_0,T,\mathbb R^m)$ and $F'w\in\mathbb H^1(t_0,T,\mathbb R^n)$) and \eqref{eq:dae} gives that 
\begin{equation}\label{eq:i1}
  \begin{split}
    \mu=&%\langle\ell,Fx(T)\rangle-\int_{t_0}^T\langle H^Tu,x\rangle dt\\
    %\langle\ell&,Fx(T)\rangle=
%\langle F'\ell,F^+Fx(T)\rangle = \langle F'w(T),F^+Fx(T)\rangle\\ 
%&=\int_{t_0}^T \langle(Fx)_t,w\rangle+\langle(F'w)_t,x\rangle dt
\langle (F^+)'F'w(t_0),Fx_0\rangle+\int_{t_0}^T \langle f,w\rangle dt\\
&+\int_{t_0}^T\langle\dfrac{d(F'w)}{dt}+C'w-H'u,x\rangle dt\,.
  \end{split}
\end{equation}
By \eqref{eq:error} and \eqref{eq:i1}:  
$\sigma(T,\ell,u)=\sup_{\mathcal R_\eta}\gamma^2+\sup_{f,x_0,x}\mu^2$. 
By Cauchy inequality $\gamma^2\le 
\int_{t_0}^{T}E\langle R\eta,\eta\rangle dt
\int_{t_0}^{T}\langle R^{-1}u,u\rangle dt$. As $E\langle R\eta,\eta\rangle=\mathrm{tr}\,(R\mathcal R_\eta)$, it follows from~\eqref{eq:eta_bounds} 
% \begin{equation*}
%   E(\int_{t_0}^{T}\langle u,\eta\rangle dt)^2\le 
% \int_{t_0}^{T}E\langle R\eta,\eta\rangle dt
% \int_{t_0}^{T}\langle R^{-1}u,u\rangle dt \le 
% \int_{t_0}^{T}\langle R^{-1}u,u\rangle dt 
% \end{equation*}
% so that
\begin{equation}
  \label{eq:worstcase_error}
  \sigma(T,\ell,u)=
\int_{t_0}^{T}\langle R^{-1}u,u\rangle dt + \sup_{f,x_0,x}\mu^2%(f,x_0,x)
\end{equation}
Assume $\hat\sigma(T,\ell)<+\infty$. Then $\sigma(T,\ell,u)<\infty$ for at least one $u$ so that $\sup_{f,x_0,x}\mu^2<\infty$. The term $\int_{t_0}^T \langle f,w\rangle dt$ in the first line of~\eqref{eq:i1} is bounded due to~\eqref{eq:G}. Thus 
\begin{equation}
  \label{eq:bound_on_w}
  \begin{split}
    \sup_{x} &\{
\langle F'w(t_0),F^+Fx(t_0)\rangle\\
&+\int_{t_0}^T\langle\dfrac{d(F'w)}{dt}+C'w-H'u,x\rangle dt\}<\infty
  \end{split}
\end{equation}
where $\sup$ is taken over all $x$ solving~\eqref{eq:dae} with $(x_0,f)\in\mathscr G$. \eqref{eq:bound_on_w} allows us to prove that there exists $z$ such that~\eqref{eq:zul} holds for the given $\ell$ and $u$. To do so we apply a general duality result\footnote{The proof of \eqref{eq:LfiG} for bounded $L$ and Banach spaces $\mathcal H_{1,2}$ can be found in \cite{Ioffe1974}} from~\cite{Zhuk2009d}: 
\begin{equation}
  \label{eq:LfiG}
  \sup_{x\in\mathscr D(L)}\{\langle\mathcal F, x\rangle, Lx\in G\}=
\inf_{b\in\mathscr D(L')}\{c(G,b),L'b=\mathcal F\}
\end{equation} 
provided (A1) $L:\mathscr D(L)\subset\mathcal H_1\to\mathcal H_2$ is a closed dense-defined linear mapping, (A2) $G\subset H_2$ is a closed bounded convex set and $\mathcal H_{1,2}$ are abstract Hilbert spaces. Define 
\begin{equation}
  \label{eq:Lx}
  \begin{split}
    &(Lx)(t)=(\dfrac{d(Fx)}{dt}-C(t)x(t),Fx(t_0)),\,x\in\mathscr D(L) \\
    &\mathscr D(L):=\{x:Fx\in\mathbb H_1(t_0,T,\mathbb R^n)\}
  \end{split}
\end{equation} 
Then $L$ is a closed dense defined linear mapping \cite{Zhuk2007} and
\begin{equation}
  \label{eq:sLz}
  \begin{split}
  (L'&b)(t)=-\dfrac{d(F'z)}{dt}-C'z,b\in\mathscr D(L'),\\
  \mathscr D&(L'):=\{b=(z,z_0):F'z\in\mathbb H_1(t_0,T,\mathbb R^m),\\
  &F'z(T)=0, z_0={F^+}'F'z(t_0)+d, F'd=0 \}
  \end{split}
\end{equation}
Setting $\mathcal F:=((F^+)'F'w(t_0),\dfrac{d(F'w)}{dt}+C'w-H'u)$ we see from~\eqref{eq:bound_on_w} that the right-hand part of~\eqref{eq:LfiG} is finite. Hence, there exists at least one $b\in\mathscr D(L')$ such that $L'b=\mathcal F$ or (using~\eqref{eq:sLz}) 
%\begin{equation}
%  \label{eq:4}
%  \inf\{c(G,(z,z_0)),
$-\dfrac{d(F'z)}{dt}-C'z=\dfrac{d(F'w)}{dt}+C'w-H'u$. %\}<+\infty
%\end{equation}
Setting $\tilde z:=(w+z)$ we obtain $
\dfrac{d(F'\tilde z)}{dt}+C'\tilde z-H'u=0$ and $F'\tilde z=F'\ell$. 
%with $\tilde z:=(w+z)$, where $b=(z,z_0)\in\mathscr D(L')$, verifying the equality in~\eqref{eq:4}. 
This proves~\eqref{eq:zul} has a solution. \\
On the contrary, let $z$ verify~\eqref{eq:zul} for the given $\ell$ and $u$. 
Then $z$ verifies conditions (W), therefore we can plug $z$ into \eqref{eq:i1} instead of $w$. 
Define\footnote{$G_1$ is a set of all $x_0,f$ such that $(x_0,f)\in\mathscr G$ and ~\eqref{eq:dae} has a solution $x$.}   $G_1:=\mathscr G\cap \mathcal R(L)$, where $
\mathcal R(L)$ is the range of the linear mapping $L$ defined by~\eqref{eq:Lx} and set $
S:=\sup_{(x_0,f)\in G_1}\langle F'z(t_0),F^+x_0\rangle+
\int_{t_0}^T\langle z,f\rangle dt$. Now, using \eqref{eq:worstcase_error} 
%integration-by-parts formula~\eqref{eq:ibp}, $E\Psi=0$ and \eqref{eq:etaSup} 
one derives easily 
\begin{equation}
  \label{eq:gerr}
%  \begin{split}
\sigma(T,\ell,u)=\int_{t_0}^{T}\langle R^{-1}u,u\rangle dt+S^2
%\end{split}
\end{equation}
Since $G_1$ is bounded, it follows that $\sigma(T,\ell,u)$ is finite. Let us prove \eqref{eq:umin}. Note that $S$ is a value of the support function of the set $G_1=\mathscr G\cap \mathcal R(L)$ on $({F'}^+F'z(t_0),z)$. To compute $S$ we note that $L$, $\mathscr G$ verify (A1), (A2) and $\mathrm{ int}\,G_1\ne\varnothing$. Thus (see \cite{Zhuk2009d}):   
\begin{equation}
  \label{eq:supGRL}
  %\sup\{\langle z,q \rangle,q\in G_1\}=%\cap \mathcal R(L)\}=
S=\min_{(z_0,z)\in N(L')}\{c(G,{F'}^+F'z(t_0)-z_0,z-v)\}
\end{equation}
and the $\min$ in~\eqref{eq:supGRL} is attained on some $(\tilde z_0,\tilde z)\in N(L')$. 
%It is easy to see that $\mathrm{ int}\,G_1\ne\varnothing$ for $L$ and $\mathscr G$ defined by~\eqref{eq:Lx} and~\eqref{eq:G} respectively. 
Recalling the definition of $L'$ (formula~\eqref{eq:sLz}) and noting that $
c^2(G,z-b)=\| Q_0^{-\frac12}({F'}^+F'(z(t_0)-v(t_0))-d)\|^2+
\int_{t_0}^T\|Q^{-\frac 12}(z-v)\|^2 dt
$ we derive \eqref{eq:umin} from \eqref{eq:gerr}-\eqref{eq:supGRL}. This completes the proof. 
% that the minimizer for \eqref{eq:umin} is also the minimizer for $\sigma(T,\ell,u)$. 
% that adjoint mapping $L^*$ is defined by the rule $$
% L^*v(t)= -\dfrac d{dt}F'v(t)-C'(t)v(t),v\in D(L^*)
% $$ where $D(L^*)=\{v\in\mathbb
% L^m_2[t_0,T]:F'v\in\mathbb W_2^n[t_0,T],F'v(T)=0\}$. Combining~\eqref{eq:dual}
% with $(*)$ and $(**)$ we obtain
% $$
% \sigma(u)=\min_v\{\int_{t_0}^T(Q^{-1}(z-v),z-v)dt\}+\int_{t_0}^T(R^{-1}u,u)dt
% $$ where $L'v=0$. This concludes the proof. 
\end{pf}
\subsection{Optimality conditions}
\label{ss:optcond}
Assume $\ell\in\mathcal L(T)$. By definition~\ref{d:1} and due to Generalized Kalman Duality (GKD) principle (see Theorem \ref{t:1}) the $\ell$-estimate $\hat u$ is a solution of the Dual Control Problem (DCP), that is the optimal control problem with cost~\eqref{eq:umin} and DAE constraint~\eqref{eq:zul} for any constant 
%\footnote{In fact, the same holds if $F(t)$ has a constant rank and $C(t)$ is measurable on $(t_0,T)$} 
$F\in\mathbb R^{m\times n}$ and continuous $t\mapsto C(t)\in\mathbb R^{m\times n}$, $t\in[t_0,T]$. If $F=I_{n\times n}$ then $\hat u=RHp$ where $p$ may be found from the following optimality conditions (Euler-Lagrange System in the Hamilton form \cite{Ioffe1974}):
\begin{equation}
    \label{eq:dae_EL}
    \begin{split}
      &\dfrac{dFp}{dt}=Cp+Q^{-1}z,\,Fp(t_0)=\tilde Q_0z(t_0),\\
      &\dfrac{dF'z}{dt}=-C'z+H'RHp,\,F'z(T)=F'\ell\,.
    \end{split}
  \end{equation}
with $\tilde Q_0 = Q_0^{-1}{F'}^+F'$. In the general case $F\in\mathbb R^{m\times n}$, let us assume that (AS) the system~\eqref{eq:dae_EL} is solvable. One can prove using direct variational method (see \cite{Ioffe1974})) that $\hat u=RHp$ solves the DCP with cost~\eqref{eq:umin} and DAE constraint~\eqref{eq:zul}. Although the assumption (AS) allows one to solve the optimal control problem with DAE constraints, %Kurina,Mermann і решту хто використовує таке припущення]. 
 it may be very restrictive for state estimation problems. 
%: forcing the solvability of~\eqref{eq:dae_EL} one constrains the minimax observability subspace $\mathcal L(T)$. 
To illustrate this, let us consider an example. Define   
\begin{equation}
  \label{eq:exmpl_FCH}
  F'=\left[
\begin{smallmatrix}
  1&&0\\
  0&&1\\
  0&&0\\
  0&&0\\
\end{smallmatrix}\right], 
 C'(t)=\left[
\begin{smallmatrix}
0&&-1\\
0&&0\\
1&&0\\
0&&-1
\end{smallmatrix}\right],
H'(t)=\left[\begin{smallmatrix}
1&&0&&0\\
0&&0&&1\\
0&&0&&0\\
0&&1&&0
\end{smallmatrix}\right]
\end{equation} and take $Q_0=Q(t)=I_{2\times 2}$, $R(t)=I_{4\times 4}$. 
In this case~\eqref{eq:dae_EL} reads as:  
% iff $\ell=(0,0)$.
\begin{equation}
  \label{eq:bvppr}
  \begin{split}
    &\dfrac{dz_1}{dt}=z_2+p_1, z_1(T)=\ell_1, -z_1=0, \\
    &\dfrac{dz_2}{dt}=p_2,z_2(T)=\ell_2, z_2+p_4=0,\\
    &\dfrac{dp_1}{dt}=p_3+z_1, p_1(t_0)=z_1(t_0),\\
    &\dfrac{dp_2}{dt}=-p_1-p_4+z_2,p_2(t_0)=z_2(t_0)\,.
  \end{split}
\end{equation}
We claim that~\eqref{eq:bvppr} has a solution iff $\ell_1=\ell_2=0$. Really, $z_1(t)\equiv0$ implies $z_1(T)=\ell_1=0$, $-z_2=p_1=p_4$ and $\frac d{dt}p_1=p_3$. According to this we rewrite~\eqref{eq:bvppr} as follows:
\begin{equation}
  \label{eq:bvpr2}
  \begin{split}
    &\dfrac{dp_1}{dt}=-p_2,p_1(t_0)=0,p_1(T)=\ell_2,\\
    &\dfrac {dp_2}{dt}=-3p_1,p_2(t_0)=0\,.
  \end{split}
\end{equation}
It is clear that~\eqref{eq:bvpr2} has a solution iff $\ell_2=0$. Thus, the assumption (AS) leads to the trivial minimax observability subspace: $\mathcal L(T)=\{0\}\times\{0\}$. However, $\mathcal L(T)=\{0\}\times\mathbb R$. To see this, take $u_3\in\mathbb L_2$, $\ell_2\in\mathbb R$ and $\ell_1=0$, and define $z_1=0$, $u_{1,2}=-z_2$, $z_2=\ell_2-\int_t^Tu_3(s)ds$. By direct substitution one checks that $z_{1,2}$ and $u_{1,2,3}$ solve~\eqref{eq:zul}. 
Therefore $\mathcal L(T)=\{0\}\times\mathbb R$ due to~\eqref{eq:LT}. We see that classical optimality condition (Euler-Lagrange system in the form~\eqref{eq:dae_EL}) may be inefficient for solving the minimax state estimation problems for DAEs.  
%we loose the information about the structure of minimax observable subspace which is contained in the DAE~\eqref{eq:zul}. 
In the next proposition we prove that optimal control problem with cost~\eqref{eq:umin} and DAE constraint~\eqref{eq:zul} has a unique solution $\hat u$, $\hat z$, provided $\ell\in\mathcal L(T)$, and we present one possible approximation of $\hat u$, $\hat z$ based on the Tikhonov regularization method \cite{Tikhonov1977}.   
%and represents its approximation. 
\begin{prop}[optimality conditions]
\label{p:3e} 
Let $\varepsilon>0$. % and $\widetilde Q_0 = Q_0^{-1}{F'}^+F'$. 
The DAE boundary-value problem 
\begin{equation}
\label{eq:daebvpe}
\begin{split}
&\frac{d(F'z)}{dt}=-C'z+H'\hat u+\hat p,\\
&\frac{d(Fp)}{dt}=Cp+\varepsilon Q^{-1}z,\,\varepsilon\hat u = RHp,\\
%Q_0^{-1}({F'}^+F'z(t_0)-d)\\
&F'z(T)+F^+Fp(T)=F'\ell,\,F'd=0\,,\\
&\frac 1\varepsilon Fp(t_0)= Q^{-1}_0({F'}^+F'z(t_0)-d)\,.
%&\varepsilon\hat u = Rp,
\end{split}
\end{equation} 
has a unique solution $\hat u_\varepsilon$, $\hat p_\varepsilon$, $\hat z_\varepsilon$, $\hat d_\varepsilon$. If $
\ell\in\mathcal L(T)$ then there exists $\hat d,\hat u$ and $\hat z$ such that 1) $\hat d_\varepsilon\to\hat d$ in $\mathbb R^n$ and $\hat u_\varepsilon\to\hat u$, $\hat z_\varepsilon\to\hat z$ in $\mathbb L_2$, 2) $\hat u$ and $\hat z$ verify \eqref{eq:zul} and 3) $\hat u$ is the $\ell$-estimate and $\ell$-error is given by
\begin{equation}
  \label{eq:mnmx_err}
  \begin{split}
    \hat\sigma(T,&\ell)=
    \sigma(T,\ell,\hat u) = \Omega(\hat u,\hat z,\hat d):=\|R^{-\frac 12}\hat u\|_{\mathbb L_2}^2\\
    &+\| Q_0^{-\frac12}({F'}^+F'\hat z(t_0)-\hat d)\|^2+\|Q^{-\frac 12}\hat z\|_{\mathbb L_2}^2\,.
  \end{split}
\end{equation}
 \end{prop}
 \begin{pf}
Define $r:=\mathrm{rang} F$ and $D=\mathrm{diag}(\lambda_1\dots\lambda_r)$ where $\lambda_i$, $i=\overline{1,r}$ are positive eigen values of $FF'$. If $r=0$ then~(\ref{eq:daebvpe}) is obviously uniquely solvable. Assume $r>0$.  
%If $n-r=0$ then the unique solvability of~\eqref{eq:daebvpe} can be proved by the standard argument (see \cite{Boichuk2004}). Let us prove that \eqref{eq:daebvpe} has a unique solution by splitting method~\cite{Eremenko1980}, provided $n-r>0$. 
%that is to split~\eqref{eq:daebvpe} into differential and algebraic parts. 
It is easy to see, applying the SVD decomposition \cite{Albert1972} to $F$, that 
$ F=U'S V$, where $UU'=I_m$, $V'V=I_{n}$ and $S=\left
  [\begin{smallmatrix}
    D^\frac12&&0_{r\times n-r}\\
    0_{m-r\times r}&&0_{m-r\times n-r}
  \end{smallmatrix}\right]$. 
Thus multiplying the first equation of~\eqref{eq:daebvpe} by $U$, the second -- by $V$, and changing variables 
one can reduce the general case to the case of DAE~\eqref{eq:daebvpe} with  
% $D=I\in\mathbb R^{r\times r}$. 
% Therefore, without loss of generality, we assume
$F=\left[
\begin{smallmatrix}
  I_{r}&&0_{r\times n-r}\\0_{m-r\times r}&&0_{m-r\times n-r}
\end{smallmatrix}
\right]$. In what follows, therefore, we can focus on this case only. 
Having in mind the above 4-block representation for $F$ 
we split the coefficients of~\eqref{eq:daebvpe} as follows: 
$C(t)=\left[
\begin{smallmatrix}
C_1&&C_2\\
C_3&&C_4  
\end{smallmatrix}
\right]$, $
Q=\left[
\begin{smallmatrix}
Q_1&&Q_2\\
Q'_2&&Q_4  
\end{smallmatrix}
\right]
$, $Q_0=\left[
\begin{smallmatrix}
Q_1^0&& Q^0_2\\
(Q^0_2)'&& Q^0_4  
\end{smallmatrix}\right]
$, $
H'RH=\left[
\begin{smallmatrix}
S_1&&S_2\\
S_2'&&S_4
\end{smallmatrix}
\right]
$, $\ell=\bigl(
\begin{smallmatrix}
  \ell_1\\
  \ell_2
\end{smallmatrix}\bigr)$, $d=\bigl(
\begin{smallmatrix}d_1\\d_2\end{smallmatrix}\bigr)$. If 1) $n-r=0$ and $m-r>0$ we set $
C_2:=0_{r\times 1}$, $C_4:=0_{m-r\times 1}$ and $S_2:=0_{r\times 1}$, $S_4:=0$; if 2) $n-r>0$ and $m-r=0$ we set $C_3:=0_{1\times r}$, $C_4:=0_{1\times n-r}$ and $Q^0_2,Q_2:=0_{m\times 1}$, $Q^0_4,Q_4:=0$; if 3) $n=m=r$ we set $C_4:=0$, $C_2:=0_{r\times 1}$, $C_3:=0_{1\times r}$ and let $S_i$, $Q^0_i,Q_i$ be defined as in 1) and 2) respectively, $i\in\{2,4\}$. 
According to this \eqref{eq:daebvpe} splits into %differential part:
%$p_1(t_0)=\varepsilon(Q^0_1z_1(t_0)-Q^0_2d_2)$, $z_1(T)+p_1(T)=\ell_1$, and 
\begin{equation*}
  \begin{split}
    &\dfrac{dp_1}{dt}=C_1p_1+C_2p_2+\varepsilon (Q_1z_1+Q_2z_2),\\
    &\dfrac{dz_1}{dt}=-C_1'z_1-C_3'z_2+p_1+\frac 1\varepsilon (S_1p_1+S_2p_2),\\
    &z_1(T)+p_1(T)=\ell_1, p_1(t_0)=\varepsilon(Q^0_1z_1(t_0)-Q^0_2d_2)\\
  \end{split}
\end{equation*}
and algebraic part: $Q^0_4d_2=(Q^0_2)'z_1(t_0)$,
\begin{equation*}
  \begin{split}
&C_3p_1+C_4p_2+\varepsilon (Q_2'z_1+Q_4z_2)=0,\\% Q^0_4d_2=(Q^0_2)'z_1(t_0),\\
&-C_2'z_1-C'_4z_2+\frac 1\varepsilon S_2'p_1+(I+\frac 1\varepsilon S_4)p_2=0.
  \end{split}
\end{equation*}
Define $\widetilde Q_4:=Q^0_1-Q^0_2(Q^0_4)^{-1}(Q^0_2)'$ and set
% $
% W(t,\varepsilon)=\varepsilon I_{n-r}+S_4+C_4'Q_4^{-1}C_4
% $, $
% M(t,\varepsilon)=W^{-1}(t,\varepsilon)
% $ and 
\begin{equation*}
  \begin{split}
&W(t,\varepsilon)=\varepsilon I_{n-r}+S_4+C_4'Q_4^{-1}C_4\,, M(t,\varepsilon)=W^{-1}(t,\varepsilon)\,,\\
&A(t)=(C_3'Q_4^{-1}C_4+S_2),\,B(t)=(C'_2-C_4'Q_4^{-1}Q_2')\,, \\
&C_\varepsilon(t)=-C_1'+C_3'Q_4^{-1}Q_2'+A(t)M(t,\varepsilon)B(t)\,,\\
% &W(t,\varepsilon)=(\varepsilon I+S_4+C_4'Q_4^{-1}C_4), M(t,\varepsilon)=
% W^{-1}(t,\varepsilon) \\
&Q_\varepsilon(t)=-\frac 1\varepsilon A(t)M(t,\varepsilon)A'(t)+I_{r}+\frac 1\varepsilon[S_1+C_3'Q_4^{-1}C_3]\,,\\
&S_\varepsilon(t)=\varepsilon (Q_1 -Q_2Q_4^{-1}Q_2'
+ B'(t)M(t,\varepsilon)B(t))\,.\\
%&M(t,\varepsilon)=(\varepsilon I+S_4+C_4'Q_4^{-1}C_4)^{-1},
  \end{split}
\end{equation*}
Solving the algebraic equations for $z_2,p_2,d_2$ we find:
\begin{equation}\label{eq:p2z2}
    \begin{split}
      &Q_4z_2= (-Q_2'-C_4MB)z_1+\frac 1\varepsilon (C_4MA'-C_3)p_1,\\
      &p_2=\varepsilon MBz_1- MA'p_1, d_2=(Q^0_4)^{-1}(Q^0_2)'z_1(t_0).
    \end{split}
  \end{equation}
Substituting~\eqref{eq:p2z2}  
%the resulting expressions 
into differential equations for $p_1,z_1$ we obtain 
\begin{equation}
\label{eq:ebvp_reduced}
  \begin{split}
    &\dfrac{dz_1}{dt}=C_\varepsilon z_1+Q_\varepsilon p_1,z_1(T)+p_1(T)=\ell_1,\\
    &\dfrac{dp_1}{dt}=-C'_\varepsilon p_1+S_\varepsilon z_1, p_1(t_0)=
\varepsilon\widetilde Q_4z_1(t_0).
  \end{split}
\end{equation}
We claim that \eqref{eq:ebvp_reduced} has a unique solution for any $\ell_1\in\mathbb R^r$ and $\varepsilon>0$. Let us prove uniqueness. Note that $Q_1 -Q_2Q_4^{-1}Q_2'>0$ as $Q(t)>0$ (see \cite{Albert1972} for details) and thus $S_\varepsilon(t)>0$ for $\varepsilon>0$ (as $B'MB\ge0$). Applying simple matrix manipulations one can prove (see, for instance, \cite{Kurina1986g}) that $Q_\varepsilon(t)\ge0$ for $\varepsilon>0$. Assume $z_1,p_1$ solve \eqref{eq:ebvp_reduced} for $\ell_1=0$. Then, integrating by parts and using \eqref{eq:ebvp_reduced} we obtain 
\begin{equation*}
  \begin{split}
    -\langle z_1(T),z_1(T)\rangle&-\langle z_1(t_0),\varepsilon\widetilde Q_4z_1(t_0)\rangle\\
    &=\int_{t_0}^T\langle Q_\varepsilon p_1,p_1\rangle+\langle S_\varepsilon z_1,z_1\rangle dt
  \end{split}
\end{equation*}
This equality is possible only if $p_1=0$, $z_1=0$ as $S_\varepsilon(t),Q_\varepsilon(t)\ge0$. As \eqref{eq:ebvp_reduced} is a Noether Boundary-Value Problem (BVP), which has a unique solution for $\ell_1=0$, it follows from the general theory of linear BVP \cite{Boichuk2004} that \eqref{eq:ebvp_reduced} has a unique solution for any $\ell_1$. Thus, we proved unique solvability of~\eqref{eq:daebvpe}. %Now assume $\ell\in\mathcal L(T)$ and let us prove 1)-2). 
Let us introduce the following definitions. Take $u\in\mathbb L_2$ and $z\in\mathbb L_2$ such that $F'z\in\mathbb H_1$, and assign to $u,z$ a number $\delta(u,z)$: 
$$  
%%%%%\delta(u,z):=\|F'z(T)-F'\ell\|^2+\int_{t_0}^T\|\dfrac{dF'z}{dt}+C'z-H'u\|^2dt 
\delta(u,z):=\|F'z(T)-F'\ell\|^2+\|\dfrac{dF'z}{dt}+C'z-H'u\|^2_{\mathbb L_2} 
$$
It was proved in \cite{Zhuk2007} that $\delta$ is convex and weakly low semi-continuous\footnote{that is $\delta(u,z)\le\lim\delta(u_n,z_n)$, provided $u_n,z_n$ converges weakly to $(u,z)$ or equally $\int_{t_0}^T\langle u_n,p\rangle+\langle z_n,q\rangle dt\to\int_{t_0}^T\langle u,p\rangle+\langle z,q\rangle dt$ for any $(p,q)$ in $\mathbb L_2\times\mathbb L_2$. 
} (w.l.s.c). Define a Tikhonov function $\mathcal T_\varepsilon(u,z,d):=\delta(u,z)+\varepsilon\Omega(u,z,d)$ with 
$\Omega(u,z,d)$ defined by~\eqref{eq:mnmx_err}. 
%:=%\|\widetilde Q_0^{-\frac12}z(t_0)-Q_0^{-\frac12}d\|^2+
%\|Q_0^{-\frac12}({F'}^+F'z(t_0)-d)\|^2+
%\|Q^{-\frac    12}z\|_{\mathbb L_2}^2+\|R^{-\frac12}u\|_{\mathbb L_2}^2dt
%$.  
% For any $u\in\mathbb L_2(t_0,T,\mathbb R^p)$ we compute using~\eqref{eq:daebvpe} 
% \begin{equation*}
%   \begin{split}
%     &T_\varepsilon(u,z,d)-T_\varepsilon(\hat u_\varepsilon,\hat z_\varepsilon,
% \hat d_\varepsilon)=\|F'z(T)-F'\ell\|^2-\|F^+F\hat p_\varepsilon(T)\|^2\\
% &+\int_{t_0}^T\|\dfrac{dF'z}{dt}+C'z-H'u\|^2+\varepsilon(\|u\|^2-\|\hat u_\varepsilon\|^2+\|z\|^2-\|\hat z_\varepsilon\|^2)dt\\
% &-\int_{t_0}^T\|\hat p_\varepsilon\|^2dt+
% %\varepsilon(\|\tilde Q_0^{-\frac12}({F'}^+F'z(t_0)-d)\|^2-\langle F\hat 
% \varepsilon(\|\tilde Q_0^{-\frac12}z(t_0)-Q_0^{-\frac12}d)\|^2-
% \langle F\hat p_\varepsilon(t_0),{F'}^+F'\hat z_\varepsilon(t_0)\rangle
%   \end{split}
% \end{equation*}
We claim (see appendix for the details) that 
\begin{equation}
  \label{eq:infTikhonovFunc}
  \mathcal T_\varepsilon(\hat u_\varepsilon,\hat z_\varepsilon,\hat d_\varepsilon)=
\inf_{u,z,d}\mathcal T_\varepsilon(u,z,d):=\mathcal T^*_\varepsilon, \quad \forall\varepsilon>0
\end{equation}
Take any $\ell\in\mathcal L(T)$. By~\eqref{eq:LT} there exists 
$u$ and $z$ such that $\inf\delta=\delta(u,z)=0$. Using \eqref{eq:infTikhonovFunc} we obtain
\begin{equation}
  \label{eq:TikhonovInequalities}
  \begin{split}
    &%0=\inf_{u,z} 
%\delta(u,z)
0\le %T(\hat u_\varepsilon,\hat z_\varepsilon,0)\le 
\mathcal T^*_\varepsilon\le \mathcal T_\varepsilon(u,z,d)
%&\le\inf_{u,z} \delta(u,z)+\varepsilon \Omega(u,z,d)+\varepsilon\\
\le \delta(\hat u_\varepsilon,\hat z_\varepsilon)+\varepsilon \Omega(u,z,d)
  \end{split}
\end{equation}
for any $d:F'd=0$. Due to \eqref{eq:TikhonovInequalities}:
$\mathcal T^*_\varepsilon=
\delta(\hat u_\varepsilon,\hat z_\varepsilon)+
\varepsilon\Omega(\hat u_\varepsilon,\hat z_\varepsilon,\hat d_\varepsilon)\le 
\delta(\hat u_\varepsilon,\hat z_\varepsilon)+\varepsilon \Omega(u,z,d)
$ so that
%$\varepsilon\Omega(\hat u_\varepsilon,\hat z_\varepsilon,\hat d_\varepsilon)\le\varepsilon \Omega(u,z,d)+\varepsilon$ implying 
\begin{equation}
  \label{eq:uezeMinimizesOmegaOverAdjDAE}
 \Omega(\hat u_\varepsilon,\hat z_\varepsilon,\hat d_\varepsilon)\le \Omega(u,z,d),\quad \delta (u,z)=0,F'd=0
\end{equation}
\eqref{eq:uezeMinimizesOmegaOverAdjDAE} proves that $\{\hat u_\varepsilon,\hat z_\varepsilon,\hat d_\varepsilon\}$ is a bounded sequence in the Hilbert space $\mathbb L_2\times\mathbb L_2\times\mathbb R^n$. Thus 
%(see \cite{Ioffe1974}) 
it contains a sub-sequence $\{\hat u_{\varepsilon_k},\hat z_{\varepsilon_k},\hat d_{\varepsilon_k}\}$ which converges weakly to some element $\hat u,\hat z,\hat d$. By~\eqref{eq:TikhonovInequalities}    
$0=\delta(u,z)\le \delta(\hat u_{\varepsilon_k},\hat z_{\varepsilon_k})\le\mathcal T_\varepsilon^*\le 
  \delta(u,z)+\varepsilon\Omega(u,z,d)$ 
so that by w.l.s.c. of $\delta$:
\begin{equation}
  \label{eq:TuepstoTopt}
\delta(\hat u,\hat z)\le
\lim_{k\to\infty}\delta(\hat u_{\varepsilon_k},\hat z_{\varepsilon_k})=\delta(u,z)  =\inf_{u,z} \delta=0
\end{equation}
%By~\eqref{eq:TuepstoTopt} $\delta(\hat u,\hat z)=0$. 
We claim that (see appendix for technical details)
\begin{equation}
  \label{eq:OmegaLowCont}
\Omega(\hat u,\hat z,\hat d)
\le\lim_{\varepsilon\to0}\Omega(\hat u_{\varepsilon},\hat z_{\varepsilon},\hat d_\varepsilon)
\end{equation} 
By \eqref{eq:OmegaLowCont} and \eqref{eq:uezeMinimizesOmegaOverAdjDAE} we get:
\begin{equation}
  \label{eq:Omegainf}
  \Omega(\hat u,\hat z,\hat d)=\Omega^*:=\inf_{\{(u,z):\delta(u,z)=0\},F'd=0}\Omega(u,z,d) 
\end{equation}
Note that $\Omega$ is strictly convex, therefore $\Omega$ has a unique minimizer $w^*$, which coincides with $(\hat u,\hat z,\hat d)$ by~\eqref{eq:Omegainf}. This proves that $w^*$ is a unique weak limiting point for the bounded sequence $\{\hat u_\varepsilon,\hat z_\varepsilon,\hat d_\varepsilon\}$. Thus, $\{\hat d_\varepsilon\}$ converges to $\hat d$ in $\mathbb R^n$ as in $\mathbb R^n$ the weak convergence is equivalent to the strong convergence. Moreover,~\eqref{eq:OmegaLowCont} and \eqref{eq:uezeMinimizesOmegaOverAdjDAE} implies $\lim\Omega(\hat u_{\varepsilon},\hat z_{\varepsilon},\hat d_\varepsilon)=\Omega(\hat u,\hat z,\hat d)$. The latter proves 1) as $\{\hat u_\varepsilon,\hat z_\varepsilon\}$ converges to $(\hat u,\hat z)$ in $\mathbb L_2$ if and only if $\{\hat u_\varepsilon,\hat z_\varepsilon\}$ converges to $(\hat u,\hat z)$ weakly and $\lim\|\hat u_\varepsilon\|^2_{\mathbb L_2}+\|\hat z_\varepsilon\|^2_{\mathbb L_2}=
\|\hat u\|^2_{\mathbb L_2}+\|\hat z\|^2_{\mathbb L_2}$ (see \cite{Ioffe1974} for details). 2) also holds as $\delta(\hat u,\hat z)=0$ by~\eqref{eq:TuepstoTopt}. 
% This and w.l.s.c. of the norm in a Hilbert space \cite{Balakrishnan1976} imply   
% \begin{equation}
%   \label{eq:ze0tohatz0}
%   \begin{split}
%     \Omega(\hat u,\hat z, \hat d)&\le 
%     %\lim_{\varepsilon\to 0}\Omega(\hat u_\varepsilon,\hat z_\varepsilon,0)
%     \lim_{\varepsilon\to 0}\|Q_0^{-\frac12}({F'}^+F'\hat z_\varepsilon(t_0)-\hat d_\varepsilon)\|^2\\
%     &+\lim_{\varepsilon\to 0}
%     \int_{t_0}^T\|Q^{-\frac    12}\hat z_\varepsilon\|^2+\|R^{-\frac12}\hat u_\varepsilon\|^2dt
%   \end{split}
% \end{equation}  
% By \eqref{eq:infTikhonovFunc} and \eqref{eq:TuepstoTopt} $
% \mathcal T^*_\varepsilon
% %\delta(\hat u_\varepsilon,\hat z_\varepsilon)+\varepsilon\Omega(\hat u_\varepsilon,\hat z_\varepsilon,\hat d_\varepsilon)
% \le\delta(\hat u,\hat z)+\varepsilon\Omega(\hat u,\hat z,\hat d)\le 
% \delta(\hat u_\varepsilon,\hat z_\varepsilon)+\varepsilon\Omega(\hat u,\hat z,\hat d)
% $ so that $\Omega(\hat u_\varepsilon,\hat z_\varepsilon,\hat d_\varepsilon) 
% \le\Omega(\hat u,\hat z,\hat d)$. This and \eqref{eq:ze0tohatz0} gives 
% \begin{equation}
%   \label{eq:strong_convergense}
%   \Omega(\hat u,\hat z, \hat d)\le \lim\Omega(\hat u_\varepsilon,\hat z_\varepsilon,\hat d_\varepsilon)\le \Omega(\hat u,\hat z, \hat d)
% \end{equation}
Let us prove 3). % and~\eqref{eq:mnmx_err}. 
% We claim that 
% $\Omega(\hat u,\hat z,\hat d)\le\Omega(u,z,d)$, provided 
% $\delta(u,z)=0$, $F'd=0$ and $\ell\in\mathcal L(T)$. To see this we note that  
% $\Omega(\hat u,\hat z,\hat d)\le\lim_\varepsilon\Omega(\hat
% u_\varepsilon,\hat z_\varepsilon,\hat d_\varepsilon)\le\Omega(u,z,d)$ due to~(\ref{eq:uezeMinimizesOmegaOverAdjDAE}) and \eqref{eq:strong_convergense}. 
Take any $u,z$ verifying~(\ref{eq:zul}) and $\tilde v,\tilde d$ defined by \eqref{eq:umin}.
% There exist $v$ verifying \eqref{eq:zul} with $u=0$ and $\ell=0$ and
% $d$ verifying $F'd=0$ such that $\min$ in \eqref{eq:umin} is
% achieved and $$ $\sigma(T,\ell,u)=\|\tilde
% Q_0^{-\frac12}(z(t_0)-v(t_0))-Q_0^{-\frac12}d\|^2
% +\int_{t_0}^T\|Q^{-\frac 12}(z-v)\|^2 dt+ \int_{t_0}^T\|
% R^{-\frac12}u\|^2dt$.  In other words, $v,d$
% minimize~(\ref{eq:umin}).
Define $\tilde z:=z-\tilde v$. 
%Recalling the definition of $\Omega$ we
%have $\sigma(T,\ell,u)=\Omega(u,\tilde z,\tilde d)$. 
Using the definition of $\tilde v$ (see \eqref{eq:umin} and notes after it) we
find that $\tilde z$ also solves~(\ref{eq:zul}). Thus, $\sigma(T,\ell,u)=
\Omega(u,\tilde z,\tilde d)\ge\Omega(\hat u,\hat z,\hat d)$ by \eqref{eq:Omegainf}.
Hence, $\hat\sigma(T,\ell)\ge\Omega(\hat u,\hat z,\hat d)$. 
% where $\min$
%is taken over all $u,z$ solving \eqref{eq:zul} and all $d:F'd=0$. 
On the other hand, we get by 1):
\begin{equation*}
  \begin{split}
    &\hat\sigma(T,\ell)\le\sigma(T,\ell,\hat u)
    \le\min_{v}\Omega(\hat u,\hat z-v,\hat d)\\
    %{\| \tilde Q_0^{-\frac12}(\hat z(t_0)-v(t_0))-Q_0^{-\frac12}\hat d\|^2\\
    %&+\int_{t_0}^T\|Q^{-\frac 12}(\hat z-v)\|^2 dt\}+\int_{t_0}^T\| R^{-\frac12}\hat u\|^2dt\\
    &\le%\mu_\varepsilon(\hat u):=
    \min_{v}\{\lim_{\varepsilon\to0}\| Q_0^{-\frac12}({F^+}'F'\hat z_\varepsilon(t_0)-{F^+}'F'v(t_0)-\hat d_\varepsilon)\|^2\\
    &+\lim_{\varepsilon\to0}\|Q^{-\frac 12}(\hat z_\varepsilon-v)\|^2\}+
    \| R^{-\frac12}\hat u\|_{\mathbb L_2}^2\\
    &=\lim_{\varepsilon\to0}\Omega(\hat u,\hat z_\varepsilon,\hat d_\varepsilon)+
    \min_{v}\{\|Q_0^{-\frac12}{F^+}'F' v(t_0)\|^2+\|Q^{-\frac 12}v\|_{\mathbb L_2}^2\}\\
    &=\Omega(\hat u,\hat z,\hat d)
%    &=\int_{t_0}^T\| R^{-\frac12}\hat u\|^2dt+\lim_{\varepsilon\to0}\| \widetilde Q_0^{-\frac12}\hat z_\varepsilon(t_0)-Q_0^{-\frac12}\hat d_\varepsilon\|^2\\
%&+\min_{v}\{\|\widetilde Q_0^{-\frac12} v(t_0)\|^2+\int_{t_0}^T\|Q^{-\frac 12}v\|^2 dt\}\\
%&+\lim_{\varepsilon\to0}\| \widetilde Q_0^{-\frac12}\hat z_\varepsilon(t_0)-Q_0^{-\frac12}\hat d_\varepsilon\|^2+
%&+\lim_{\varepsilon\to0}\int_{t_0}^T\|Q^{-\frac 12}\hat z_\varepsilon\|^2 dt=\Omega(\hat u,\hat z,\hat d)%=\min_{u,z,d}\Omega(u,z,d)
  \end{split}
\end{equation*}
where we obtained the 4th line noting that $\|Q^{-\frac 12}(\hat z_\varepsilon-v)\|_{\mathbb L_2}^2 = 
    \|Q^{-\frac 12}\hat z_\varepsilon\|_{\mathbb L_2}^2 +
    \|Q^{-\frac 12}v\|_{\mathbb L_2}^2 -
    2\int_{t_0}^T\langle Q^{-1}\hat z_\varepsilon,v\rangle dt$ and 
   \begin{equation*}
  \begin{split}
&\int_{t_0}^T\langle Q^{-1}\hat z_\varepsilon,v\rangle dt
    =\frac 1\varepsilon \int_{t_0}^T\langle\dfrac{dF\hat p_\varepsilon}{dt}-C\hat p_\varepsilon,v\rangle dt\\
&=-\langle Q_0^{-1}({F^+}'F'\hat z_\varepsilon(t_0)-\hat d_\varepsilon),{F^+}'F'v(t_0)\rangle
\end{split}
\end{equation*}
where the latter equality follows from~(\ref{eq:ibp}), definition of $v$ (see notes after~\eqref{eq:umin}) and~(\ref{eq:daebvpe}). 
% , it follows that 
% \begin{equation*}
%   \begin{split}
% \mu_\varepsilon(\hat u)=&
% \int_{t_0}^T\| R^{-\frac12}\hat u\|^2dt+\lim_{\varepsilon\to0}\| \widetilde Q_0^{-\frac12}\hat z_\varepsilon(t_0)-Q_0^{-\frac12}\hat d_\varepsilon\|^2\\
% &+\min_{v}\{\|\widetilde Q_0^{-\frac12} v(t_0)\|^2+\int_{t_0}^T\|Q^{-\frac 12}v\|^2 dt\}\\
% %&+\lim_{\varepsilon\to0}\| \widetilde Q_0^{-\frac12}\hat z_\varepsilon(t_0)-Q_0^{-\frac12}\hat d_\varepsilon\|^2+
% &+\lim_{\varepsilon\to0}\int_{t_0}^T\|Q^{-\frac 12}\hat z_\varepsilon\|^2 dt=\Omega(\hat u,\hat z,\hat d)%=\min_{u,z,d}\Omega(u,z,d)
%   \end{split}
% \end{equation*}
Thus $\hat\sigma(T,\ell)\le\sigma(T,\ell,\hat u)\le\mu_\varepsilon(\hat u)=\Omega(\hat u,\hat z,\hat d)$ implying~\eqref{eq:mnmx_err}. \eqref{eq:mnmx_err} 
%$\hat\sigma(T,\ell)=\sigma(T,\ell,\hat u)=\Omega(\hat u,\hat z,\hat d)$. 
implies, in turn, $\hat u$ is the $\ell$-estimate by definition. This completes the proof.  
 \end{pf}
We will refer $\hat u_\varepsilon$ as a sub-optimal $\ell$-estimate. 
Let us represent $\hat u_\varepsilon(y)$ in the form of the minimax filter.  
%Using SVD decomposition (see the first part of the proof of Proposition~\ref{p:3e}) one can always transform the general case to the case 
% $F=\bigl[
% \begin{smallmatrix}
%   I_{r}&&0_{r\times n-r}\\
% 0_{m-r,r}&&0_{m-r\times m-r}
% \end{smallmatrix}
% \bigr]$.
Recalling definitions of $M,A,B$ introduced at the beginning of the proof of Proposition~\ref{p:3e}, and splittings for $\ell,Q,R,H,C$ we define $\Phi(t,\varepsilon)= \left[
\begin{smallmatrix}
  K(t,\varepsilon)\\
M(t,\varepsilon)[\varepsilon B(t)-A'(t)K(t,\varepsilon)]
\end{smallmatrix}\right]$, where $K(t,\varepsilon)$ solves $$
\frac {dK}{dt}=-KC_\varepsilon-C'_\varepsilon K-
   KQ_\varepsilon K+S_\varepsilon,\,K(0)=\varepsilon \tilde Q_4
$$ Define $z_1(T)=(I_r+K(T,\varepsilon))^{-1}\ell_1$ and let $z_1$ solve 
\begin{equation}
  \label{eq:zK}
  \begin{split}
    &  \frac{dz_1}{dt}=C_\varepsilon z_1+Q_\varepsilon K(t,\varepsilon)z_1,\\
 \end{split}
\end{equation}
Define $\hat x_\varepsilon(t_0)=0$  
%$z_1(T)=(I_r+K(T,\varepsilon))^{-1}\ell_1$ 
%and $K(0,\varepsilon)=
%\varepsilon \tilde Q_4$ 
and let $\hat x_\varepsilon$ solve the following linear stochastic differential equation: 
\begin{equation}\label{eq:xKz}
  \begin{split}
   &d\hat x_\varepsilon = [-C'_\varepsilon-K(\varepsilon,t)Q_\varepsilon]
   \hat x_\varepsilon dt + \Phi H'R dy,\\%\hat x_\varepsilon(0)=0\\
   % &\frac{dz_1}{dt}=C_\varepsilon z_1+Q_\varepsilon K(t,\varepsilon)z_1,\\%,z_1(T)=(I+K_\varepsilon(T))^{-1}\ell_1\\
%   &\frac {dK}{dt}=-KC_\varepsilon-C'_\varepsilon K-
%   KQ_\varepsilon K+S_\varepsilon\,K(0)=\varepsilon \tilde Q_4
%   &\hat x(T):=\lim_{\varepsilon\to0}\varepsilon^{-1}(I+K(T,\varepsilon))^{-1}E\hat x^\varepsilon(T)
    \end{split}
  \end{equation}
\begin{cor}%[sub-optimal $\ell$-estimate]
\label{c:eme}
Assume $F=\bigl[
\begin{smallmatrix}
  I_{r}&&0_{r\times n-r}\\
0_{m-r,r}&&0_{m-r\times m-r}
\end{smallmatrix}
\bigr]$. % and $n-\mathrm{rank }F>0$. 
%Define $\hat x(T):=\lim_{\varepsilon\to0}\varepsilon^{-1}(I+K(T,\varepsilon))^{-1}\hat x^\varepsilon(T)$. Then 1) 
%1) $\hat x(T)$ is well defined, 
%2) %the $\ell$-estimate can be represented as 
%$E(\int_{t_0}^T\langle\hat u,y\rangle dt)=\langle \ell_1,\hat x(T)\rangle$ for any $\ell\in\mathcal L(T)$,  
Then %sub-optimal $\ell$-estimate can be represented as 
$$
\hat u_\varepsilon(y)=\int_{t_0}^T\langle\hat u_\varepsilon{},y\rangle{dt}=\langle\varepsilon^{-1}(I_r+K(T,\varepsilon))^{-1}\hat x_\varepsilon(T),\ell_1\rangle
$$ The sub-optimal $\ell$-error is given by  
$$
\hat\sigma^\varepsilon(T,\ell):=\frac1\varepsilon\bigr[
\langle z_1(T),z_1(T)
%(I+K(T,\varepsilon))^{-1}\ell_1,
%K(T,\varepsilon)(I+K(T,\varepsilon))^{-1}
%\ell_1
\rangle
-\int_{t_0}^T\|\Phi(t,\varepsilon)z_1\|^2dt\bigl]
$$ 
\end{cor}
%  \begin{rem}
%    Note, that the first equation in \eqref{eq:xKz} is a linear stochastic equation as $y$ depends on $\eta$ which is a realization of a random process verifying \eqref{eq:eta_bounds}. \eqref{eq:eta_bounds} means that the correlation function $E\eta(t)\eta^T(t)$ is bounded. This imply that 
% % $\mathrm 
% % E\int_{t_0}^T\langle\eta(t),\eta{(t)}\rangle dt\le 1$
% the integral (in the sense described in \cite{Loev1968})
% $\int_{t_0}^t\eta(s)ds$ exists for any $t\in[t_0,T]$ so that the integral of $y$ and its linear deterministic transformations exists. Therefore, the random process $\hat x_\varepsilon(t)$ may be defined by the formula $$
% \hat x_\varepsilon(t) = 
% \int_{t_0}^t \Upsilon(t,s)\Phi H'R y(s)ds
% $$ where $\Upsilon_t(t,s)=[-C'(\varepsilon,t)-K(\varepsilon{},t)Q(\varepsilon,t)]\Upsilon$, $\Upsilon(s,s)=I$. Roughly speaking, the realization $y$ belongs to $\mathbb L_2(t_0,T,\mathbb R^p)$ with probability 1, so that $\hat x_\varepsilon(t)$ is defined with probability 1.    
%  \end{rem}
 \begin{pf}
Assume $\hat p_\varepsilon$ solves \eqref{eq:daebvpe}. We split $p_\varepsilon=\bigl(
\begin{smallmatrix}
  p_1\\p_2
\end{smallmatrix}\bigr)$ where $p_1$ solves \eqref{eq:ebvp_reduced} and $p_2$ is defined by \eqref{eq:p2z2}. It can be checked by direct calculation that $p_1(t)=K(t,\varepsilon)z_1(t)$ where 
%$p_1$ solves \eqref{eq:ebvp} and 
$z_1$ is defined by~\eqref{eq:zK}. Using this and \eqref{eq:p2z2} we deduce $\hat p_\varepsilon
= \Phi(t,\varepsilon)z_1$. \eqref{eq:daebvpe} implies $\hat u_\varepsilon=\frac 1\varepsilon RH\hat p_\varepsilon$. Finally, using the obtained representations for $\hat p_\varepsilon$, $\hat u_\varepsilon$ and~\eqref{eq:xKz} we obtain integrating by parts that
\begin{equation*}
%\label{eq:xe_approximates_hatu}
  \begin{split}
    \hat u_\varepsilon&(y)=\int_{t_0}^T\langle y,\frac 1\varepsilon RH\hat p_\varepsilon\rangle dt = \int_{t_0}^T\langle \frac 1\varepsilon \Phi' H'R y, z_1\rangle dt\\
&=\langle \ell_1,\varepsilon^{-1}(I_r+K(T,\varepsilon))^{-1}\hat x_\varepsilon(T)\rangle\to
\int_{t_0}^T\langle \hat u,y\rangle dt=\hat u(y)
  \end{split}
\end{equation*}
By~\eqref{eq:uezeMinimizesOmegaOverAdjDAE}-\eqref{eq:OmegaLowCont} 
%As it was proved in Proposition~\ref{p:3e}, 
$\Omega(\hat u_\varepsilon,\hat z_\varepsilon,\hat d_\varepsilon)\to\hat\sigma(T,\ell)$. It is easy to compute using~\eqref{eq:daebvpe} that $\Omega(\hat u_\varepsilon,\hat z_\varepsilon,\hat d_\varepsilon)=\varepsilon^{-1}(\langle \ell,F\hat p_\varepsilon\rangle+\|\hat p_\varepsilon\|_{\mathbb L_2}^2)$. To conclude it is sufficient to substitute $\hat p_\varepsilon= \Phi(t,\varepsilon)z_1$ into the latter formula. % for $\Omega$.
\end{pf}
\textbf{Example.} In order to demonstrate main benefits of Proposition~\ref{p:3e} we will apply it to the example presented above: assume that the bounding set, state equation and observation operator are defined by~\eqref{eq:exmpl_FCH}. Note that ${F'}^+F'=I_4$. Thus $F'd=0$ implies $d=0$. According to Theorem~\ref{t:1} the exact $\ell$-estimate $\hat u = \bigl(
\begin{smallmatrix}
  \hat u_1\\
  \hat u_2\\
  \hat u_3
\end{smallmatrix}\bigr)$ may be obtained minimizing  $
\sigma(T,\ell,u)=\sum_{i=1}^2\|z_i(t_0)\|^2+\|z_i\|_{\mathbb L_2}^2+
\sum_{j=1}^3\|u_j\|^2_{\mathbb L_2}
$ over solutions of the DAE 
\begin{equation}
  \label{eq:w12}
  \begin{split}
    &\dfrac {dz_1}{dt}-z_2-u_1=0, z_1(T)=\ell_1, -z_1=0\\
    &\dfrac{dz_2}{dt}-u_3=0,z_2(T)=\ell_2,-z_2-u_2=0
  \end{split}
\end{equation}
Assume $\ell_1=0$ so that $\ell=\bigl(
\begin{smallmatrix}
  \ell_1\\
  \ell_2
\end{smallmatrix}\bigr)\in\mathcal L(T)$. If $\hat u_{1,2}$ solves~\eqref{eq:w12} then $\hat u_{1,2}=-z_2$. Hence, $\hat u_3$ may be found minimizing $
\sigma(T,\ell,u)=\|z_2(t_0)\|^2+\int_{t_0}^T 3z_2^2+u_3^2dt
$ over $\dfrac{dz_2}{dt}=u_3$, $z_2(T)=\ell_2$. The optimality condition takes the following form: $\hat u_3=p$, $\frac{dz_2}{dt}=p$, $z_2(T)=\ell_2$, $\frac{dp}{dt}=3z_2$, $p(t_0)=z_2(t_0)$. Let us represent the estimate in the form of the minimax filter. Introducing $k$ as a solution of the Riccati equation $\frac{dk}{dt}=k^2-3$, $k(0)=1$ we find that $\hat u_2=kz_2$ where $z_2$ solves the following Cauchy problem: $\frac {dz_2}{dt} = kz_2$, $z_2(T)=\ell_2$. Let $\hat x$ be a solution to $\frac{d\hat x}{dt}=-k\hat x-y_1-y_2+ky_3$, $\hat x(t_0)=0$. Then it is easy to see that $\hat u(y)=\int_{t_0}^T\langle \hat u, y\rangle dt = \ell_2\hat x(T)$, where $y=\bigl(
\begin{smallmatrix}
  y_1\\y_2\\y_3
\end{smallmatrix}\bigr)$ denotes a realization of the random process representing observed data. Let us compute the sub-optimal $\ell$-estimate. \eqref{eq:xKz} reads as 
\begin{equation*}
  \begin{split}
    \dfrac d{dt}&
    \bigl(
    \begin{smallmatrix}
      x_1\\x_2
    \end{smallmatrix}
    \bigr)
    =\left[\begin{smallmatrix}
       -(1+\frac 1\varepsilon )k_1&&-(1+\frac 1\varepsilon )k_2\\
       -1-(1+\frac 1\varepsilon )k_2&&-(1+\frac 1\varepsilon )k_4
      \end{smallmatrix}\right]
    \bigl(
    \begin{smallmatrix}
      x_1\\x_2
    \end{smallmatrix}
    \bigr)\\
&+\left[\begin{smallmatrix}
        k_1&&k_2&&1&&0\\
        k_2&&k_4&&0&&\frac{-\varepsilon}{1+\varepsilon}
      \end{smallmatrix}\right]
    \left[\begin{smallmatrix}
        1&&0&&0\\
        0&&0&&1\\
        0&&0&&0\\
        0&&1&&0
      \end{smallmatrix}\right]
    \bigl(
    \begin{smallmatrix}
      y_1\\y_2\\y_3
    \end{smallmatrix}
    \bigr),\bigl(
    \begin{smallmatrix}
      x_1(t_0)\\x_2(t_0)
    \end{smallmatrix}
    \bigr)=\bigl(
    \begin{smallmatrix}
      0\\0
    \end{smallmatrix}
\bigr)
\\
\dfrac d{dt}&\left[\begin{smallmatrix}
        k_1&&k_2\\
        k_2&&k_4
      \end{smallmatrix}\right]=
    \left[\begin{smallmatrix}
        0&&0\\
        -1&&0
      \end{smallmatrix}\right]
    \left[\begin{smallmatrix}
        k_1&&k_2\\
        k_2&&k_4
      \end{smallmatrix}\right]+
\left[\begin{smallmatrix}
        k_1&&k_2\\
        k_2&&k_4
      \end{smallmatrix}\right]
\left[\begin{smallmatrix}
        0&&-1\\
        0&&0
      \end{smallmatrix}\right]\\
&\left[\begin{smallmatrix}
        1+\varepsilon&&0\\
        0&&\varepsilon+\frac\varepsilon{1+\varepsilon}
      \end{smallmatrix}\right]-
(1+\frac 1\varepsilon)
\left[\begin{smallmatrix}
        k_1&&k_2\\
        k_2&&k_4
      \end{smallmatrix}\right]^2
,\left[\begin{smallmatrix}
        k_1(t_0)&&k_2(t_0)\\
        k_2(t_0)&&k_4(t_0)
      \end{smallmatrix}\right]=\left[\begin{smallmatrix}
        1&&0\\
        0&&1
      \end{smallmatrix}\right].
  \end{split}
\end{equation*}
Define $\bigl(
    \begin{smallmatrix}
      \tilde \ell_1\\\tilde \ell_2
    \end{smallmatrix}
    \bigr)=\left[\begin{smallmatrix}
        1+k_1(T)&&k_2(T)\\
        k_2(T)&&1+k_4(T)
      \end{smallmatrix}
\right]^{-1}\bigl(
    \begin{smallmatrix}
      \ell_1\\\ell_2
    \end{smallmatrix}
    \bigr)$. 
Due to Corollary~\ref{c:eme} the sub-optimal $\ell$-estimate may be represented as $\hat u_\varepsilon(y)=\int_{t_0}^T\langle \hat u_\varepsilon, y\rangle dt = \frac 1\varepsilon
\langle \bigl(
    \begin{smallmatrix}
      \tilde \ell_1\\\tilde \ell_2
    \end{smallmatrix}
    \bigr),\bigl(
    \begin{smallmatrix}
      x_1(T)\\x_2(T)
    \end{smallmatrix}
    \bigr)\rangle$. If $\ell_1=0$ then the sub-optimal $\ell$-error is given by 
    $\hat\sigma^\varepsilon(T,\ell)=\frac 1\varepsilon
\langle\left[\begin{smallmatrix}
        k_1(T)&&k_2(T)\\
        k_2(T)&&k_4(T)
      \end{smallmatrix}
\right]
\bigl(
    \begin{smallmatrix}
      \tilde \ell_1\\\tilde \ell_2
    \end{smallmatrix}
    \bigr),
\bigl(
    \begin{smallmatrix}
      \tilde \ell_1\\\tilde \ell_2
    \end{smallmatrix}
    \bigr)\rangle$. 
 \begin{figure}[t] 
   \centering
\includegraphics[width=.4\textheight]{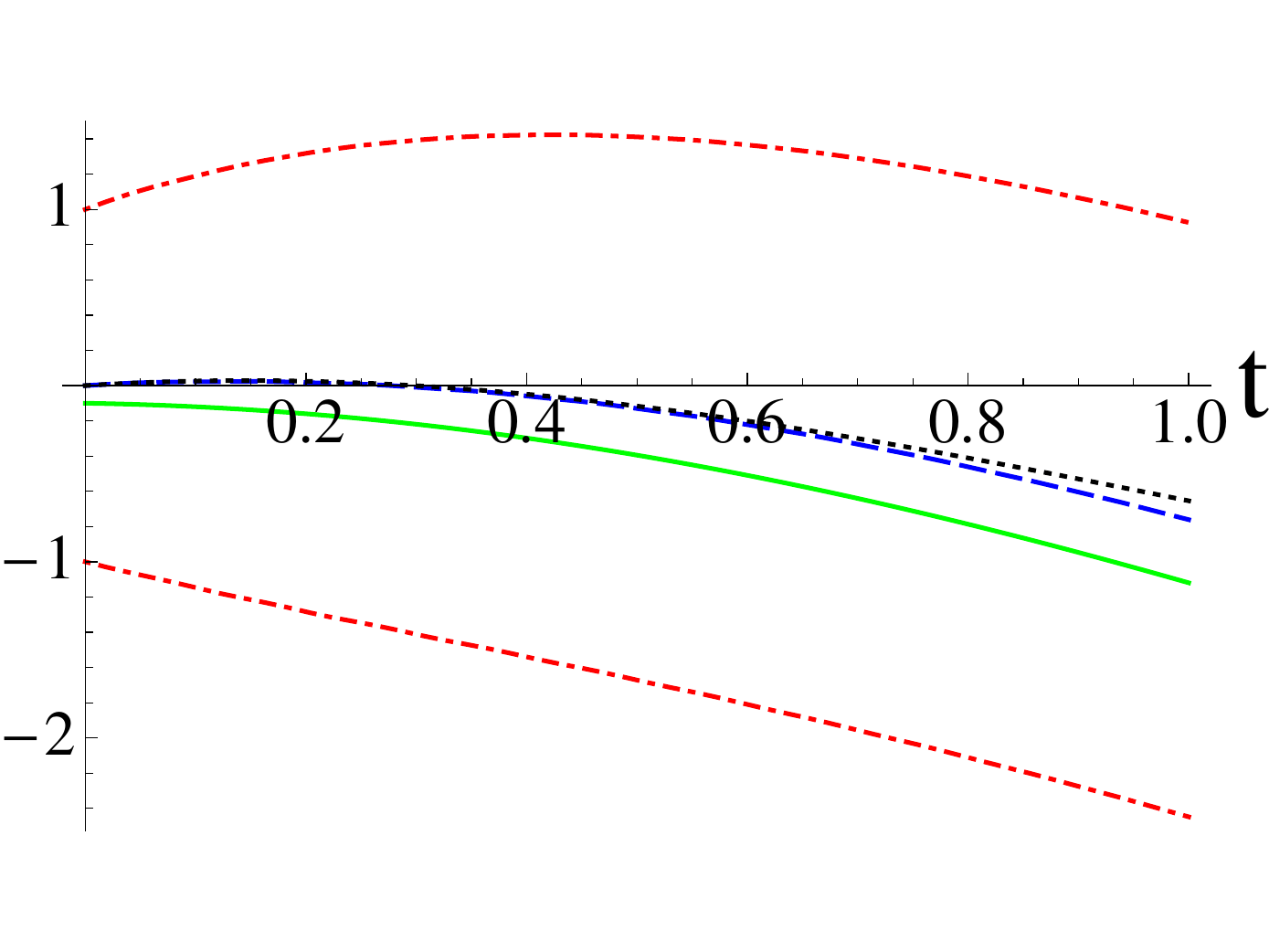}   
   \caption{
Optimal $\ell$-estimate $\hat u(y)$ (Dashed), sub-optimal $\ell$-estimate $\hat u_\varepsilon(y)$ (Dotted) and error $\hat\sigma^\varepsilon(t,\ell)$ (DotDashed) and simulated $x_2(t)$ (Solid), $\ell=(0,1)^T$, $t\in[0,1]$, $\varepsilon=\exp(-30)$.
}
   \label{fig:2}
 \end{figure}
Take $t_0=0$, $T=1$ and assume that $F$ and $C$ are defined by~\eqref{eq:exmpl_FCH}. In the corresponding DAE $x_{3,4}$ are free components. For simulations we choose $x_3=\cos(t)$ and $x_4=\sin(t)$, $x_1(0)=0.1$, $x_2(0)=-0.1$, $f_1=f_2=0$. In order to generate artificial observations $y$ we take $\eta(t)=\bigl(
\begin{smallmatrix}
  -0.1\\-0.2\\0.3
\end{smallmatrix}
\bigr)$. In Figure~\ref{fig:2} the optimal $\ell$-estimate, sub-optimal $\ell$-estimate and sub-optimal $\ell$-error are presented, provided $\ell_1=0$, $\ell_2=1$. As $\mathcal L(t)\equiv\{0\}\times\mathbb R$ we see that $x_1$ is not observable in the minimax sense. This can be explained as follows. The derivative $x_3$ of $x_1$ may be any element of $\mathbb L_2$. As we apply integration by parts formula in order to compute $\sigma(T,\ell,u)$ (see~\eqref{eq:ibp}), the expression for $\sigma(T,\ell,u)$ contains $\int_{t_0}^Tx_3z_1dt$. Thus, $\sigma(T,\ell,u)<+\infty$ implies $z_1(t)\equiv 0$ for any $t\in[t_0,T]$, in particular $z_1(T)=\ell_1=0$. In this case $\hat u\in\mathbb L_2$ and $\hat u_\varepsilon\to\hat u$ in $\mathbb L_2$. If $\ell_1\ne0$ then the only candidate for $\hat u_1$ is the impulse control $\delta(T-t)\ell_1$ switching $z_1$ from $0$ to $\ell_1$ at time-instant $T$. However, in this case 
%$\hat u_1=\delta(T-t)\ell_1\not\in\mathbb L_2$ implying that the corresponding $\ell$-error is infinite
%\footnote{
% if $\hat u_1=\delta(T-t)\ell_1$ then $\|\hat u_\varepsilon\|_{\mathbb L_2}\to\infty$, $\varepsilon\to0$ so that 
% $\Omega(\hat u_\varepsilon,\hat z_\varepsilon,\hat d_\varepsilon)\to\infty$. Thus, by \eqref{eq:uezeMinimizesOmegaOverAdjDAE} $\ell$-error
%(sub-optimal $\ell$-error depends on the $\mathbb L_2$-norm of $\hat u_\varepsilon$ (see~\eqref{eq:umin})) 
%goes to infinity. Note that 
the numerical sub-optimal $\ell$-error computed by the algorithm of Corollary~\ref{c:eme} increases: $\hat\sigma^\varepsilon(1,\ell)\approx 3\times 10^6$, provided $\ell_1=1$, $\ell_2=0$ and $\varepsilon=\exp(-30)$. 
Let us illustrate this. Assume $ C'(t)=\left[
\begin{smallmatrix}
0&&-1\\
0&&0\\
c(t)&&0\\
0&&-1
\end{smallmatrix}\right]$, where $c(t)=0$ for $t_1\le t\le T$ and $c(t)>0$ for $t_0\le t<t_1$.  Then~\eqref{eq:zul} is solvable for any $\ell$: the solution is given by $u_1:=\frac {dz_1}{dt}-z_2$, $u_2=-z_2$, $z_1(t)=0$, $t_0\le t\le t_1$ and $z_1(t)=\frac{t-t_1}{T-t_1}\ell$ for $t_1\le t\le T$. As $\hat u_1=\frac {dz_1}{dt}-z_2$, it follows that $\|\hat u_1\|_{\mathbb L_2}$ goes to infinity if $t_1\to T$. We stress that the limiting case $t_1=T$ with $c(t)\equiv 1$ corresponds to $C'(t)$ which is being considered in our example (see~\eqref{eq:exmpl_FCH}).  

\section{Conclusion}
 \label{sec:Conclusion}
The paper presents one way to generalize the minimax state estimation approach for linear time-varying DAE~\eqref{eq:dae}. The only restriction we impose here is that $F$ does not depend on time. But our approach can be generalized to the class of time-varying $F(t)$ with constant (or piece-wise constant) rank by means of Lyapunov-Floke theorem. The main idea behind the generalization is the Generalized Kalman Duality (GKD) principle. GKD allows to formulate a Dual Control Problem (DCP) which gives an expression for the Worst-Case Error (WCE). Due to GKD the WCE is finite if and only if 
Also GKD gives necessary and sufficient conditions (in the form of the minimax observable subspace) for the WCE finiteness. In order to compute $\ell$-estimate one needs to solve the DCP, that is a linear-quadratic control problem with DAE constraints. Application of the classical optimality conditions (Euler-Lagrange equations) imposes additional constraints onto the minimax observability subspace $\mathcal L(T)$. To avoid this we apply a Tikhonov regularization approach allowing to construct sub-optimal solutions of DCP or, sub-optimal estimates. If $\ell\in\mathcal L(T)$ then the sequence of sub-optimal estimates converges to the $\ell$-estimate which belongs to $\mathbb L_2$. Otherwise sub-optimal estimates weakly converge to the linear combination of delta-functions. The $\mathbb L_2$-norms of the sub-optimal $\ell$-estimates grow infinitely in this case. 
%In perspective, it is important to define a generic procedure for choosing the regularization parameter $\varepsilon$ depending on the time discretization step for DAE. General remarks on this topic are available in \cite{Tikhonov1977}. Another key point is to find an efficient  description of the minimax observability subspace $\mathcal L(T)$ for the case of time-varying $C(t)$. In particular, it is important to have easy to check sufficient conditions for $\mathcal L(T)=\mathbb R^m$. The case $C(t)\equiv C$ may be addressed by the method proposed in \cite{Gantmacher1960}. 

\textbf{Appendix.} 
\appendix
Let us prove~\eqref{eq:infTikhonovFunc}. Integrating by parts (formulae~\eqref{eq:ibp}) one finds
\begin{equation*}
  \begin{split}
\int_{t_0}^T&\langle \dfrac{dF\hat p_\varepsilon}{dt}-C\hat p_\varepsilon,z\rangle dt= -\int_{t_0}^T\langle\hat p_\varepsilon,\dfrac{dF' z}{dt}+C'z\rangle dt\\
&+\langle F\hat p_\varepsilon(T),{F'}^+F'z(T)\rangle-
\langle F\hat p_\varepsilon(t_0),{F'}^+F'z(t_0)\rangle
  \end{split}
\end{equation*}
In particular
\begin{equation*}
  \begin{split}
    &\int_{t_0}^T\langle \dfrac{dF\hat p_\varepsilon}{dt}-Cp_\varepsilon,\hat z_\varepsilon\rangle dt=\langle F'\hat z_\varepsilon(T)-F'\ell,F'\hat z_\varepsilon(T)\rangle\\
    &-\int_{t_0}^T\langle\hat p_\varepsilon,\hat p_\varepsilon+H'\hat u_\varepsilon\rangle dt\\
&-\langle Q_0^{-1}({F'}^+F'\hat z_\varepsilon(t_0)-d_\varepsilon),
{F'}^+F'\hat z_\varepsilon(t_0)-\hat d_\varepsilon\rangle
  \end{split}
\end{equation*}
Having this in mind it is straightforward to check that
\begin{equation*}
  \begin{split}
    &\frac\varepsilon2\int_{t_0}^T(\|R^{-\frac 12}u\|^2-\|R^{-\frac 12}\hat u_\varepsilon\|^2+\|Q^{-\frac 12}z\|^2-\|Q^{-\frac 12}\hat z_\varepsilon\|^2dt)\\
&\ge\int_{t_0}^T\langle\varepsilon R\hat u_\varepsilon,u-\hat u_\varepsilon\rangle
+\langle\varepsilon Q\hat z_\varepsilon, z-\hat z_\varepsilon\rangle dt=
\int_{t_0}^T\langle H\hat p_\varepsilon,u-\hat u_\varepsilon\rangle dt\\
&+\int_{t_0}^T\langle \dfrac{dF\hat p_\varepsilon}{dt}-Cp_\varepsilon,z\rangle dt -
\int_{t_0}^T\langle \dfrac{dF\hat p_\varepsilon}{dt}-Cp_\varepsilon,\hat z_\varepsilon\rangle dt=\\
&\int_{t_0}^T(\|\hat p_\varepsilon\|^2-\langle\hat p_\varepsilon,\dfrac{dF' z}{dt}+C'z-H'u\rangle) dt\\
&+\langle F\hat p_\varepsilon(T),{F'}^+F'z(T)\rangle-\langle F\hat p_\varepsilon(t_0),{F'}^+F'z(t_0)\rangle\\
&-\langle F\hat p_\varepsilon(T),{F'}^+F'\hat z_\varepsilon(T)\rangle+
\langle F\hat p_\varepsilon(t_0),{F'}^+F'\hat z_\varepsilon(t_0)\rangle
  \end{split}
\end{equation*}
where we have applied the sub-gradient inequality \cite{Rockafellar1970} to pass from the first line to the second line. Using this inequality, the definition of $\mathcal T_\varepsilon$ and~\eqref{eq:daebvpe} it is straightforward to check that 
\begin{equation*}
  \begin{split}
    &\mathcal T_\varepsilon(u,z,d)-\mathcal T_\varepsilon(\hat u_\varepsilon,\hat z_\varepsilon,\hat d_\varepsilon)
    \ge\|F'z(T)-F'\ell+F^+F\hat p_\varepsilon(T)\|^2\\
    &+\|Q_0^{-\frac12}({F'}^+F'z(t_0)-{F'}^+F'\hat z_\varepsilon(t_0)+\hat d_\varepsilon-d)\|^2\\
    &+\int_{t_0}^T\|\dfrac{dF'z}{dt}+C'z-H'u-\hat p_\varepsilon\|^2dt\ge 0
  \end{split}
\end{equation*}
%The latter estimate proves~\eqref{eq:infTikhonovFunc}. 
Let us prove~\eqref{eq:OmegaLowCont}. We proved that $\hat u_{\varepsilon},\hat z_{\varepsilon}$ converges weakly to $\hat u,\hat z$ and $\{\hat d_\varepsilon\}\to\hat d$ in $\mathbb R^n$. As the norm in $\mathbb L_2$ is weakly low 
semi-continuous, it follows that
\begin{equation*}
  \begin{split}
    \Omega(\hat u_{\varepsilon},&\hat z_{\varepsilon}, \hat d_{\varepsilon})-\|Q_0^{-\frac12}({F'}^+F'z_{\varepsilon}(t_0)-d_{\varepsilon})\|^2\\
&\le \Omega(\hat u,\hat z,\hat d)-\|Q_0^{-\frac12}({F'}^+F'\hat z(t_0)-\hat d)\|^2
  \end{split}
\end{equation*}
Therefore it is sufficient to show that $F'z_\varepsilon(t_0)\to F'\hat z(t_0)\text{ in }\mathbb R^n$. Noting that $F'q(t_0)=F'q(T)-\int_{t_0}^T \dfrac{dF'q}{dt}(t)dt$ for any $q\in\mathbb H_1$ we write
\begin{equation}\label{eq:Fze(T)toFz(T)}
  \begin{split}
    &|\langle F'z_\varepsilon(t_0)-F'\hat z(t_0),v\rangle|\le 
    \|F'z_\varepsilon(T)-F'\hat z(T)\|\times\\
    &\times\|v\|+|\int_{t_0}^T\langle\dfrac{dF'\hat z_\varepsilon}{dt}-\dfrac{dF'\hat z}{dt},v\rangle dt|
  \end{split}
\end{equation}
(\ref{eq:TuepstoTopt}) implies $\|F'z_\varepsilon(T)-F'\hat z(T)\|\to0$ and 
%the sequence of numbers $\{\delta(\hat u_\varepsilon,\hat z_\varepsilon)\}$ is bounded. 
\begin{equation}
  \label{eq:dFzCzHuBounded}
  \int_{t_0}^T\|\{\dfrac{dF'\hat z_\varepsilon}{dt}+C'\hat z_\varepsilon-H'\hat u_\varepsilon\}\|^2dt<+\infty,
\,\forall\varepsilon>0
\end{equation}
%$\{\dfrac{dF'\hat z_\varepsilon}{dt}+C'\hat z_\varepsilon-H'\hat u_\varepsilon\}$ is bounded. 
%contains weakly convergent sub-sequence. 
As $\hat z_\varepsilon$ and $\hat u_\varepsilon$ converge weakly, it follows that $\lim\{C'\hat z_\varepsilon-H'\hat u_\varepsilon\}=C'\hat z-H'\hat u$. This and (\ref{eq:dFzCzHuBounded}) implies $\{\dfrac{dF'\hat z_\varepsilon}{dt}\}$ is bounded. Therefore, the weak convergence of $\hat z_\varepsilon$ gives: 
\begin{equation}
  \label{eq:dFze2dFz}
  \lim_{\varepsilon\to0}\int_{t_0}^T\langle\dfrac{dF'\hat z_\varepsilon}{dt},v\rangle dt=
\int_{t_0}^T\langle\dfrac{dF'\hat z}{dt},v\rangle dt,\,v\in\mathbb L_2
\end{equation} 
\eqref{eq:dFze2dFz} implies $\langle F'z_\varepsilon(t_0)-F'\hat z(t_0),v\rangle$ in~\eqref{eq:Fze(T)toFz(T)} converges to zero for any $v\in\mathbb R^n$ implying $F'z_\varepsilon(t_0)\to F'\hat z(t_0)$. 
%\label{sec:appendix}
%$Q(t,\varepsilon)\ge 0$, $S(t,\varepsilon)\ge 0$. %Morrison-Woodbery ???

\bibliographystyle{apalike}
\bibliography{refs/refs,refs/myrefs}

\begin{thebibliography}{}

\bibitem[Albert, 1972]{Albert1972}
Albert, A. (1972).
\newblock {\em Regression and the Moor-Penrose pseudoinverse}.
\newblock Acad. press, N.Y.

\bibitem[Boichuk and Samoilenko, 2004]{Boichuk2004}
Boichuk, A. and Samoilenko, A. (2004).
\newblock {\em Generalized Inverse Operators and Fredholm Boundary-Value
  Problems}.
\newblock VSP, Utrecht.

\bibitem[Campbell, 1987]{Campbell1987}
Campbell, S. (1987).
\newblock A general form for solvable linear time varying singular systems of
  differential equations.
\newblock {\em SIAM J. Math. Anal.}, 18(4).

\bibitem[Chernousko, 1994]{Chernousko1994}
Chernousko, F.~L. (1994).
\newblock {\em State Estimation for Dynamic Systems .}
\newblock Boca Raton, FL: CRC.

\bibitem[Darouach et~al., 1997]{Darouach1997}
Darouach, M., Boutayeb, M., and Zasadzinski, M. (1997).
\newblock Kalman filtering for continuous descriptor systems.
\newblock In {\em ACC}, pages 2108--2112, Albuquerque. AACC.

\bibitem[Gantmacher, 1960]{Gantmacher1960}
Gantmacher, F. (1960).
\newblock {\em The theory of matrices}.
\newblock Chelsea Publish.Comp., N.-Y.

\bibitem[Gerdin et~al., 2007]{Gerdin2007}
Gerdin, M., Sch{\"o}n, T.~B., Glad, T., Gustafsson, F., and Ljung, L. (2007).
\newblock On parameter and state estimation for linear differential-algebraic
  equations.
\newblock {\em Automatica J. IFAC}, 43(3):416--425.

\bibitem[Ioffe and Tikhomirov, 1974]{Ioffe1974}
Ioffe, A. and Tikhomirov, V. (1974).
\newblock {\em Theory of extremal problems}.
\newblock North-Holland, Amsterdam.

\bibitem[Kurina, 1986]{Kurina1986g}
Kurina, G.~A. (1986).
\newblock Linear {H}amiltonian systems that are not solved with respect to the
  derivative.
\newblock {\em Differentsialnye Uravneniya}, 22(2):193--198, 363.

\bibitem[Kurzhanski and V{\'a}lyi, 1997]{Kurzhanski1997}
Kurzhanski, A. and V{\'a}lyi, I. (1997).
\newblock {\em Ellipsoidal calculus for estimation and control}.
\newblock Systems \& Control: Foundations \& Applications. Birkh\"auser Boston
  Inc., Boston, MA.

\bibitem[Milanese and Tempo, 1985]{Tempo1985}
Milanese, M. and Tempo, R. (1985).
\newblock Optimal algorithms theory for robust estimation and prediction.
\newblock {\em IEEE Trans. Automat. Control}, 30(8):730--738.

\bibitem[Nakonechny, 1978]{Nakonechnii1978}
Nakonechny, A. (1978).
\newblock A minimax estimate for functionals of the solutions of operator
  equations.
\newblock {\em Arch. Math. (Brno)}, 14(1):55--59.

\bibitem[Rockafellar, 1970]{Rockafellar1970}
Rockafellar, R. (1970).
\newblock {\em Convex analysis}.
\newblock Princeton.

\bibitem[Tikhonov and Arsenin, 1977]{Tikhonov1977}
Tikhonov, A. and Arsenin, V. (1977).
\newblock {\em Solutions of ill posed problems}.
\newblock Wiley, New York.

\bibitem[Zhuk, 2007]{Zhuk2007}
Zhuk, S. (2007).
\newblock Closedness and normal solvability of an operator generated by a
  degenerate linear differential equation with variable coefficients.
\newblock {\em Nonlinear Oscil.}, 10(4):464--480.

\bibitem[Zhuk, 2009]{Zhuk2009d}
Zhuk, S. (2009).
\newblock Estimation of the states of a dynamical system described by linear
  equations with unknown parameters.
\newblock {\em Ukrainian Math. J.}, 61(2):214--235.

\bibitem[Zhuk, 2010]{Zhuk2009c}
Zhuk, S. (2010).
\newblock Minimax state estimation for linear discrete-time
  differential-algebraic equations.
\newblock {\em Automatica J. IFAC}, 46(11):1785--1789.

\end{thebibliography}
%\bibliography{refs/refs,refs/myrefs}
\end{document}